\documentclass[preprint,review,12pt]{elsarticle}
\usepackage[cmex10]{amsmath}
\usepackage{amssymb}
\usepackage{graphics}
\usepackage{graphicx}
\usepackage{verbatim}
\usepackage{epsfig,float}
\usepackage{algorithm}
\usepackage{algorithmic}
\usepackage{natbib}
\usepackage{url}
\usepackage{slashbox}
\usepackage{array}
\usepackage{amsfonts}
\usepackage{layout}
\usepackage{nomencl}
\usepackage{etoolbox}
\usepackage{booktabs}
\usepackage{multirow}
\usepackage{colortbl}
\usepackage{fancyhdr}
\usepackage{romannum}
\usepackage{subcaption}
\usepackage{pbox}
\usepackage[mathscr]{euscript}

\makeatletter
\newcommand*{\rom}[1]{\expandafter\@slowromancap\romannumeral #1@}
\makeatother

\renewcommand\nomgroup[1]{%
  \item[\bfseries
  \ifstrequal{#1}{A}{Sets and Indices}{%
  \ifstrequal{#1}{B}{Input Parameters}{%
  \ifstrequal{#1}{C}{Primary Decision Variables}{%
  \ifstrequal{#1}{D}{Auxiliary Variables}
  }}}%
]}

\journal{Transportation Research Part C}
\begin{document}
\begin{frontmatter}
\title{Stochastic and Simulation-based Models for Setting Flow Rates in Collaborative Trajectory Options Program (CTOP)}

\author[isu]{G.~Zhu\corref{cor1}}
\ead{gzhu@iastate.edu}
\author[isu]{P.~Wei}
\ead{pwei@iastate.edu}
\author[metron]{R.~Hoffman}
\ead{hoffman@metronaviation.com}
\author[metron]{B.~Hackney}
\ead{hackney@metronaviation.com}

\cortext[cor1]{Corresponding author}

\address[isu]{Aerospace Engineering Department, Iowa State University, Ames, IA, 50011, U.S.A.}
\address[metron]{Metron Aviation, Herndon, VA, 20171, U.S.A.}

\begin{abstract}
As a new tool in the NextGen portfolio, the Collaborative Trajectory Options Programs (CTOP) combines multiple features from its forerunners including Ground Delay Program (GDP), Airspace Flow Program (AFP) and reroutes, and can manage multiple Flow Constrained Areas (FCAs) with a single program. A key research question in CTOP is how to set traffic flow rates under traffic demand  and airspace capacity uncertainties. In this paper, we first investigate existing CTOP related stochastic optimization models and point out their roles in CTOP flow rate planning, and their advantages and disadvantages in terms of model flexibility, performance, practicality and Collaborative Decision Making (CDM) software compatibility, etc. CTOP FCA rate planning problem has been split into two steps: traffic flow rate optimization given demand estimation, and flow rate adaptation when flight rerouting is considered. Second, we discuss in detail a class of models called FCA-PCA (Potentially Constrained Area) models, which are extended from GDP models to solve the first step of the problem, and were considered promising as they are designed to be consistent with current CTOP software implementation. We will reveal one inherent shortcoming suffered by FCA-PCA models and show that how this deficiency can be addressed by the PCA model family. We will talk about the problems that prevent stochastic programming being optimal in the second step of the problem. Third, we discuss the applicability of simulation-based optimization, combine it with stochastic programming based heuristics and test the resulting new model on a realistic use case. The results are very encouraging. The models and discussions in this work are not only useful in more effectively implementing and analyzing CTOP programs, but are also valuable for the general multiple constrained airspace resources optimization problem and the design of future air traffic flow program.
\end{abstract}

\begin{keyword}
Air Traffic Flow Management \sep Stochastic Programming Models \sep Simulation-based Optimization \sep Collaborative Trajectory Options Program
\end{keyword}
\end{frontmatter}

\makenomenclature
\mbox{}

\nomenclature[A0]{$\mathcal{F}$}{Set of FCAs}
\nomenclature[A0a]{$\mathcal{P}$}{Set of PCAs}
\nomenclature[A1]{$\mathcal{R}$}{Set of resources, including FCAs and PCAs, $\mathcal{R} = \mathcal{F \bigcup P}$}
\nomenclature[A2]{$\mathcal{C}$}{Set of ordered pairs of resources. $(r,r^\prime)\in \mathcal{C}$ iff $r$ is connected to $r^\prime$ in directed graph}
\nomenclature[A4]{$\mathcal{T}$}{Set of time periods, $t=1,\cdots,"|\mathcal{T}"|$}
\nomenclature[A5]{$\mathcal{Q}$}{Set of scenarios, $q=1,\cdots, "|\mathcal{Q}"|$}
\nomenclature[Aaa]{$\mathcal{S}$}{Set of stages, $s=1, \cdots, "|\mathcal{S}"|$}
\nomenclature[Ac]{$\mathcal{B}$}{Set of branches in the scenario tree, $b = 1,\dots, "|\mathcal{B}"|$}
\nomenclature[Ad]{$\mathscr{P}$}{Set of paths, $\rho=1,\cdots, "|\mathscr{P}"|$}

\nomenclature[B1]{$\Delta^{r,r^\prime}$}{Number of time periods to travel from resource $r$ to $r^\prime$. Defined for all pairs $(r,r^\prime) \in \mathcal{C}$ }
\nomenclature[B2]{$p_q$}{Probability that scenario $q$ occurs}
\nomenclature[B3]{$f_{t}^{r,r^\prime}$}{Fraction of flights from resource $r$ directed to resource $r^\prime$ in time period $t$}
\nomenclature[B4]{$M_{t,q}^{r}$}{Maximum capacity of PCA $r\in \mathcal{P}$ in time period $t$ under scenario $q$}
\nomenclature[B5]{$c_g, c_a$}{Cost for taking one unit of ground delay/air delay}
\nomenclature[B6]{$t_s$}{Time period in which stage $s$ begins}
\nomenclature[B7]{$S_{s,t}^r$}{Number of flights originally scheduled to depart in stage $s$ and arrive in FCA $r \in \mathcal{F}$ (direct demand) in time period $t$}
\nomenclature[B8]{$N_b$}{Number of scenarios corresponding to branch $b$}
\nomenclature[BA]{$S^r_{t,\rho}$}{Scheduled direct demand for PCA $r$ in time period $t$ from flights with the same path $\rho$ }

\nomenclature[C0]{$P^{r}_{t}$}{Planned acceptance rate at FCA $r \in \mathcal{F}$ in time period $t$}
\nomenclature[C1]{$P^{r}_{t,q}$}{Planned acceptance rate at FCA $r \in \mathcal{F}$ in time period $t$ under scenario $q$ }
\nomenclature[C2]{$X^{r,q}_{s,t,t^\prime}$}{Number of flights originally scheduled to depart in stage $s$ and arrive in FCA $r \in \mathcal{F}$ in time period $t$, rescheduled to arrive in time period $t^\prime$ under scenario $q$}
\nomenclature[C3]{$P^r_{t,\rho}$}{Planned direct demand at PCA $r$ in time period $t$ from flights with the same path $\rho$ }

\nomenclature[D0]{$L_{t,q}^{r}$}{Number of flights that actually cross PCA $r \in \mathcal{P}$ in time period $t$ under scenario $q$ }		
\nomenclature[D1]{$A^{r}_{t,q}$}{Number of flights taking air delay before crossing PCA $r \in \mathcal{P}$ in time period $t$ under scenario $q$}
\nomenclature[D2]{$\text{UpFCA}^{r}_{t}$}{Number of flights from all upstream FCAs arrived at resource $r \in \mathcal{R}$ in time period $t$}
\nomenclature[D3]{$\text{UpFCA}^{r}_{t,q}$}{Number of flights from all upstream FCAs arrived at resource $r \in \mathcal{R}$ in time period $t$ under scenario $q$}
\nomenclature[D4]{$\text{UpPCA}^{r}_{t,q}$}{Number of flights from all upstream PCAs arrived at PCA $r \in \mathcal{P}$ in time period $t$ under scenario $q$}
\nomenclature[D6]{$A^{k,q}_{t,\rho}$}{Number of flights with the same path $\rho$ taking air delay before PCA $k$ in time period $t$ under scenario $q$}
\nomenclature[D7]{$L^{k,q}_{t,\rho}$}{Number of flights with the same path $\rho$ that actually cross PCA $k$ in time period $t$ under scenario $q$ }
\nomenclature[D8]{$\text{UpPCA}^{k,q}_{t,\rho}$}{Number of flights in path $\rho$ arriving at PCA $k$ in time period $t$ from the upstream PCA in the same path under scenario $q$}

\printnomenclature
\section{Introduction}
Collaborative Trajectory Options Programs (CTOP) is a new Traffic Management Initiative (TMI) which is used by air traffic managers to balance demand with available airspace capacity. Compared with its predecessors Ground Delay Program (GDP), Airspace Flow Program (AFP) and reroutes, the introduction of CTOP brings two major benefits: the ability to handle multiple Flow Constrained Areas (FCAs) with a single program, and the flexibility it gives to airspace users to express their conditional preference over different route choices by allowing them to submit a set of desired reroute options (Trajectory Options Set or TOS) \cite{FAAWebsiteCTOP}.

An important research question in TMI optimization has been determining airport/FCA planned acceptance rates, which is the maximum number of aircraft to be accepted in each time period. Because planned acceptance rates have to be set hours in advance so that flights can absorb delays on the ground or reroute to avoid the congested airspace, when the capacity reduction is caused by weather activity, decision making has to deal with the uncertain nature the weather forecast. Most of the literature on TMI optimization has focused on single airport GDP planning, and the dominant decision making under uncertainty approach has been stochastic programming, which minimizes the expected cost under different weather scenarios. Since the 1990s, the Federal Aviation Administration (FAA) has made significant changes in their air traffic flow management, moving from a centralized system to one called Collaborative Decision Making (CDM). Many decision support tools for air traffic managers and airline personnel are developed under this CDM paradigm \cite{vossen2012air}. Representative work differs in two aspects: the  degree to which traffic managers can modify or revise flights' controlled departure times and the compatibility with current CDM software. In \cite{ball2003stochastic}, \cite{richetta1993solving}, the ground delay decisions are made at the beginning of the planning horizon and the two models are CDM-compatible \cite{kotnyek2006equitable}. In \cite{richetta1994dynamic}, the ground delay decisions are made at the flights'  original scheduled departure stage to make use of the updated weather information. In \cite{mukherjee2007dynamic}, ground delay decisions consider the fact this flight may be further ground delayed later on (“plan to replan”). Both \cite{richetta1994dynamic} and \cite{mukherjee2007dynamic} are not compatible with current CDM software, since these two models assume that the directly control of individual and group of flights is possible.

CTOP traffic flow rates planning is substantially more challenging than GDP or AFP rate planning due to two reasons. First, because there are now multiple FCAs, the locations of FCAs can be in parallel or in serial and the rate of one FCA may affect the traffic volume to adjacent or downstream FCAs. Therefore the rates for CTOP FCAs need to be determined in an integrated way. Second, in general demand estimation and capacity information are both needed to determine the optimal FCA rates \cite{ZhuPARTMI}. However, in CTOP the air traffic demand for the constrained regions is uncertain, since one flight now may have more than one route option.

This work is an endeavor to optimize CTOP traffic flow rates under uncertainties in an effectively and efficient way. The main contributions of this paper are as follows:

\begin{enumerate}
  \item We conduct a literature review and categorize the existing CTOP related stochastic optimization models based on five criteria:
  \begin{enumerate}
    \item How dynamic or flexible the model is in terms of assigning delays (and reroutes)
    \item Whether the model will only assign delay or both reroute and delay
    \item Whether FCA, the flow control mechanism, and Potentially Constrained Areas (PCA), the actual constrained airspace region, are both explicitly used in model formulation
    \item Aggregate level of decision variables
    \item CDM-CTOP software compatibility
  \end{enumerate}
  The main results are listed in Table \ref{CTOPTable}. This table not only elucidates the nuances of different CTOP models, but also can guide researchers to develop new models for other classical and future TMIs, which is valuable for air traffic flow management research.

  \begin{table}[h!]
  \centering
  \resizebox{1.0\textwidth}{!}{%
  \begin{tabular}{cccccc}
    \toprule
    &\multicolumn{2}{c}{Aggregate Models} & \multicolumn{2}{c}{Disaggregate Models}\\
    \cmidrule(lr){2-3} \cmidrule(lr){4-5}
               & FCA-PCA Rate Planning                 & PCA Rate Planning           & \multicolumn{2}{c}{Flight Level Planning for PCA Network}\\
    \cmidrule(lr){2-2} \cmidrule(lr){3-3} \cmidrule(lr){4-5}
               & Multi-commodity Flow Approximation& Multi-commodity Flow    & Lagrangian& Lagrangian-Eulerian\\
    \cmidrule(lr){2-2} \cmidrule(lr){3-3} \cmidrule(lr){4-4} \cmidrule(lr){5-5}
    Static       & ESOM\cite{HoffmanESOM}                    &   \cite{ZhuAggregate} \S \Romannum{4}      &  \cite{ZhuCentralized} \S \Romannum{3}A &\cite{ZhuCentralized}\S\Romannum{3}B\\
    Semi-dynamic & {\cellcolor[gray]{0.8}}Semi-dynamic ESOM  &   \cite{ZhuAggregate} \S \Romannum{5}      &  \cite{ZhuCentralized} \S \Romannum{4}A &\cite{ZhuCentralized}\S\Romannum{4}B\\
    Dynamic      & D-ESOM\cite{MetronMarch2018Report}        &   \cite{ZhuAggregate} \S \Romannum{6}      &  \cite{ZhuCentralized} \S \Romannum{5}A &\cite{ZhuCentralized}\S\Romannum{5}B\\
  \bottomrule
  \end{tabular}}
  \caption{Comparison of Stochastic Models for CTOP}\label{CTOPTable}
\end{table}
  \item In the CDM paradigm, CTOP rate planning problem has been split into two steps: traffic flow rate optimization given demand estimation, and flow rate adaptation when flight rerouting exercised by CTOP slot allocation algorithm is also considered.
  By investigating two instances of FCA-PCA type models, we point out that this type of models, which are designed to solve the first step of the problem and designed to be most consistent with current CDM-CTOP software, are flawed due to traffic flow approximation. We show that how this deficiency can be addressed by PCA models, and how should FCAs be created, and filters and rates be correctly set. We discuss that why second step of the problem is more difficult to handle in the stochastic programming framework due to demand shift and conservative FCA rate policy issues.

  \item Simulation-based optimization will directly evaluate the proposed FCA rates by running CTOP slot allocation and flow simulation algorithms, and can avoid the difficulties faced by stochastic optimization in the second step of the problem. However, it now faces curse of dimensionality and has to report a solution given rather limited time, e.g. 5 minutes. We talk about some guidelines and tradeoffs in applying simulation-based optimization to CTOP FCA rate planning problem.
\end{enumerate}

This paper is organized as follows: in section \ref{Preliminaries}, we introduce some important concepts which are used throughout this paper. From section \ref{CTOPReview} to section \ref{Simulation_based_Optimization}, we talk in detail items 1 to 3 listed above. In section \ref{Experiment}, we test and compare the several stochastic and simulation-based models on a realistic use case. In section \ref{Conclusions}, we summarize the major findings and contributions of our work.

\section{Preliminary Concepts}\label{Preliminaries}
In this section, we talk about seven key components of various models described in the following sections.
\begin{figure}
    \centering
    \begin{minipage}[t]{.42\textwidth}
        \centering
        \includegraphics[width=\linewidth]{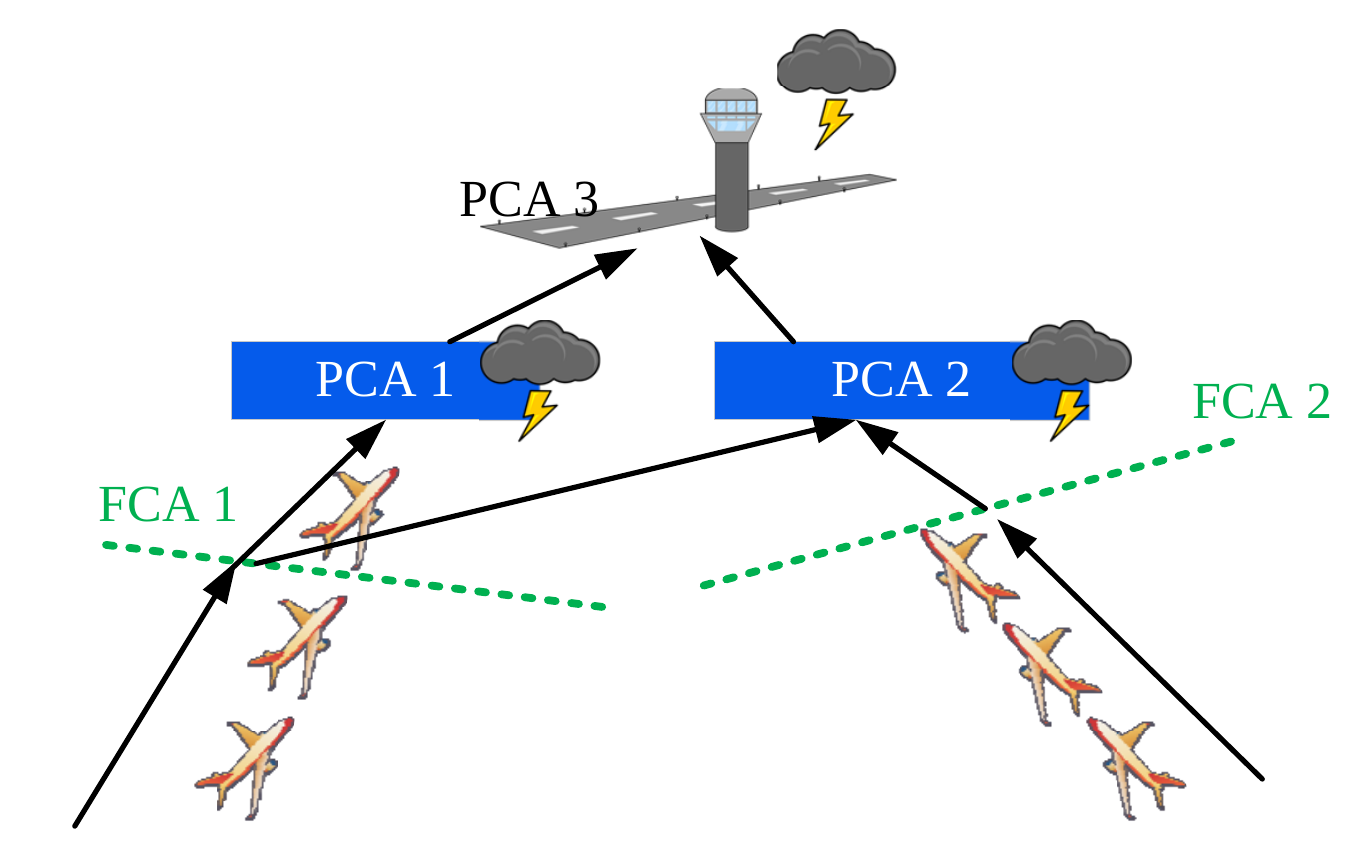}
            \caption{FCAs and PCAs}\label{FCAsPCAs}
    \end{minipage}%
    \begin{minipage}[t]{.58\textwidth}
        \centering
        \includegraphics[width=\linewidth]{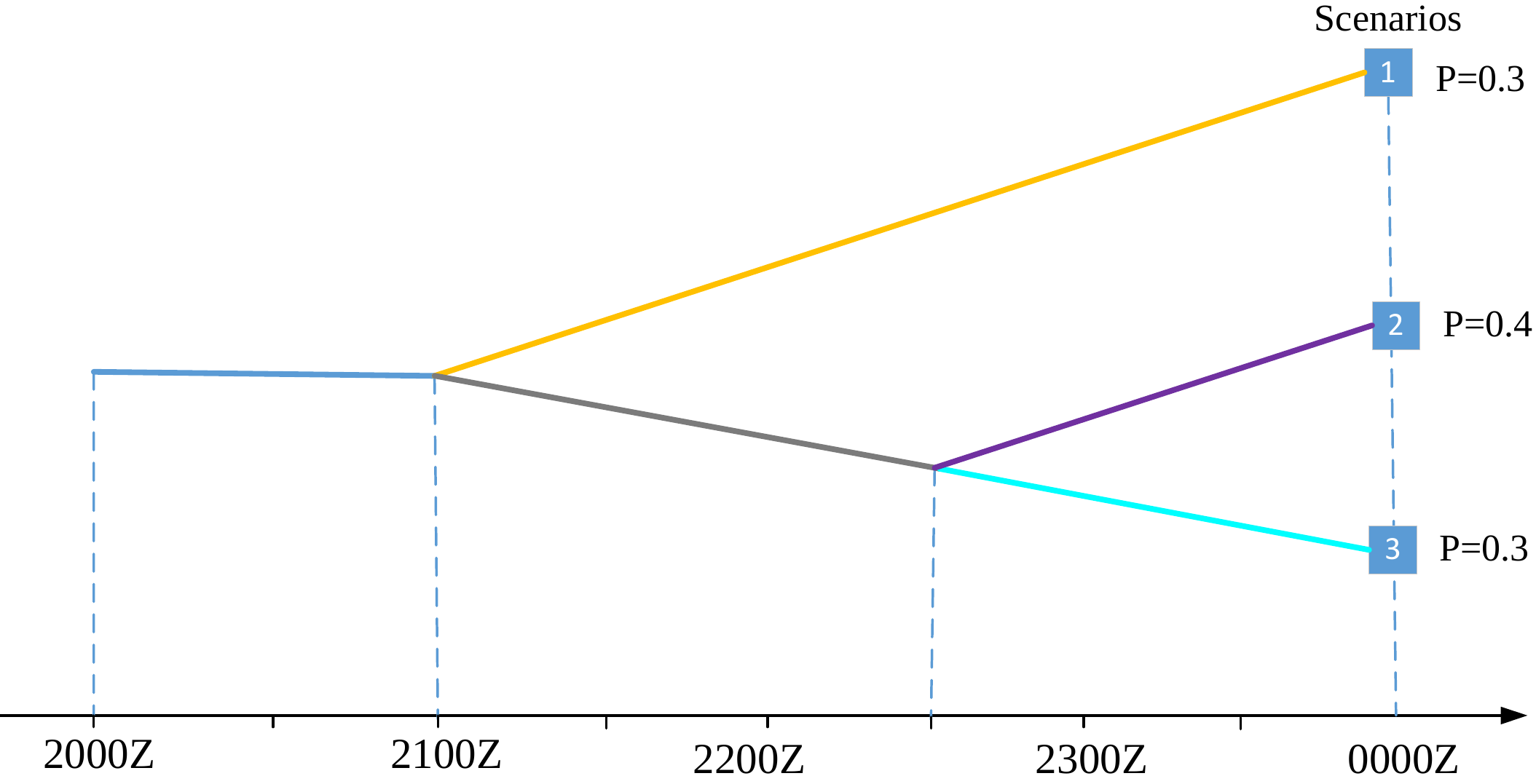}
        \caption{Scenario Tree}\label{ScenarioTreeExperiment}
    \end{minipage}
\end{figure}
\subsection{Potential Constrained Area and Capacity Scenarios}
Potential Constrained Areas (PCAs) are the airspace regions in which air traffic demand may exceed capacity. Though an imbalance resulting from pure demand increase is possible, we usually focus on the case in which adverse weather causes capacity reduction. In the stochastic programming framework, future capacity realization is represented by a finite set of scenarios arranged in a scenario tree, and each scenario is associated with a probability. A branch point in a scenario tree corresponds to the time when we acquire new information and similar scenarios evolve into distinct scenarios. For example, there are three scenarios in the scenario tree we use in this paper, shown in Figure \ref{ScenarioTreeExperiment}. Scenario 1 to 3 correspond to optimistic, average and pessimistic weather forecasts, respectively.
\subsection{Flow Constrained Area and Planned Acceptance Rate}
Flow Constrained Area (FCA) was first introduced in AFP to model constrained airspace resources. In the CTOP setting, different from a PCA which coincides with a physically constrained area and whose future capacity is stochastic, a FCA is an artificial line or region in the airspace and serves like a valve to control traffic flows into a region. Figure \ref{FCAsPCAs} depicts both FCAs and PCAs in a traffic management setting. The goal of this paper is to determine the best FCA rates on behalf of by air traffic manager. In CTOP, we explicitly distinguish the control mechanism (FCA with one set of acceptance rates) with the source of the problem (PCA with multiple sets of possible capacity values), which facilitates the multiple constrained resources mathematical modeling.

\subsection{FCA-PCA Network and Path}
A related concept is the FCA-PCA network, which is a directed graph that links the FCAs and PCAs and models the potential movement of traffic between them. This concept gives structure to our problem and allows us to use various network flow optimization techniques. General multi-resource air traffic management is by nature a multi-commodity problem, since flights will traverse different congested regions and reach different destinations. The concept of path is used to establishes a \emph{commodity}, which is the sequence of PCA nodes in the FCA-PCA network that flights traverse. For example, in Figure \ref{PCA_Network}, $\text{PCA1}\rightarrow \text{PCA}\textunderscore\text{EWR}$ is one path and $\text{PCA1}\rightarrow \text{PCA}\textunderscore\text{Exit}$ is another path. Flights in these two paths share the capacity resource of PCA1.

\subsection{Traffic Flow Rates}
To clarify the definitions of different ``rates'', we have FCA (conditional) planned acceptance rates or FCA rates for short for FCA-PCA Models; we have PCA (conditional) planned acceptance rates or PCA rates for PCA models; when we talk about traffic flow control in CTOP in general, we vaguely call FCA and PCA rates the CTOP traffic flow rates, or CTOP rates.

\section{Review and Classification of Existing CTOP Stochastic Models}\label{CTOPReview}
In this section, we discuss in detail about these twelve CTOP related stochastic models listed in Table \ref{CTOPTable} from five perspectives.
\subsection{Dynamic (Multistage) Vs. Semi-dynamic (Multistage) Vs. Static (Two-stage) Model}
The first criterion that differentiates these models is how flexibly the flights' controlled departure times (and reroutes) can be modified or revised. Similar to single airport GDP planning, if the ground delay (and reroute) are determined only using scenario information but not the branch point information, and usually assigned at the beginning of the planning horizon, the model is classified as a static model; if the ground delay (and reroute) decisions are made at a flight's original scheduled departure stage or other predetermined time, and uses the latest weather information, the model is semi-dynamic; if the ground delay (and reroutes) can be revised multiple times and each decision takes these future revisions into account, then the model is dynamic. The essence of the difference is how much information about the structure of the scenario tree and flight schedule is exploited by a model. Both the semi-dynamic models and dynamic models are multistage stochastic models, and static models are two-stage stochastic models. The dynamic models are more flexible than semi-dynamic models, which in turn are more flexible than static models.  The more flexible a model is, the better system delay performance it can achieve, but the less predictable the flight schedule will be.

\subsection{Delay Assignment Model Vs. Delay and Route Assignment Model}
The second criterion is whether the model will only assign delay or both reroute and delay. For models in columns 1 and 2 of Table \ref{CTOPTable}, we assume the routes are known, and we need to manage traffic demand through the congested regions by dynamically assigning delays to flights to minimize system delay costs,. For models in columns 3 and 4, we have the freedom to choose routes for flights and we aim to show the best total system performance we can potentially achieve in terms of total reroute and delay costs.

One of the two key features of CTOP is that it will assign not only the ground delay but also the reroute. Therefore, models in columns 1 and 2 only solve the first step of the CTOP planning problem, and have to be paired with TOS allocation algorithm (section \ref{slot_allocation_algorithm})\cite{CTOPAlgorithm}. Models in columns 3 and 4, on the other hand, solve the ground delay and reroute assignment problem at the same time in a centralized way.

\subsection{FCA-PCA Model Vs. PCA Model}
The third criterion is whether we use both FCA and PCA concepts in the model. In the current CDM-CTOP software implementation, only the concept of FCA is considered. Therefore only if a model can generate FCA planned acceptance rates, can it be used by the current software. FCA-PCA models can directly optimize and output FCA rates (section \ref{FCA_PCA_Model_section}), in this sense they are more attractive than PCA models.

Without using the artificial FCA concept at the beginning, it turns out that the mathematical formulations of PCA models are neater and the ``physical pictures'' are clearer. After solving PCA models, it is possible to obtain FCA rates by post-processing optimal PCAs rates. More details will be discussed in section \ref{ESOM_discussion} and \ref{PCA_Optimality}.
\subsection{Aggregate Model Vs. Disaggregate Model}
The fourth criterion is the aggregate level of decision variables. The CTOP rate planning problem is naturally a multi-commodity problem, since flights will traverse different congested airspace regions and reach different destinations. We can choose to use the pre-calculated flow split ratio to approximate the traffic flows between resources (column 1), or choose to explicitly deal with a multi-commodity flow problem (column 2). In the latter case, we group flights by the PCAs they traverse. The former one is more aggregate than the latter one. In multistage semi-dynamic and dynamic models, we need to further group flights by departure stages or en route times, and models will become even less aggregate.

Models in columns 3 and 4 are essentially at a flight-by-flight level and therefore called disaggregate models. They have to be at flight level since these models also determine the reroute assignment and the composition of TOSs is different for different flights.
\subsection{CDM-CTOP Compatibility}
The fifth criterion is the compatibility with current CDM-CTOP software. It is important to first define what is CDM-compatibility. In the literature, one definition of CDM-compatibility is that the model should be able to accommodate the FAA and airline operations including slots compression, intra-airline cancellation and substitution \cite{mukherjee2007dynamic}. In this loose sense, all twelve models can be made to be consistent with CDM philosophy. Here we adopt a stricter definition: the model should be compatible with the current CDM-CTOP software implementation.

All semi-dynamic and dynamic models are not compatible with the current CDM-CTOP software, since direct control individual and group of flights are not allowable and the current software does not support conditional delay decisions, same as in the GDP case. On a side note, a practical reason is that compared with static models is that a weather forecast scenario tree with relatively accurate branch points is not easily obtainable hitherto for multiple constrained en route resources, which is an important line of research in the aviation weather community \cite{liu2008scenario, steiner2010translation, clarke2012determining, LL2015AirspaceFlowRate}. Disaggregate models, even though they do not directly optimize FCA rates, can provide guidelines on setting FCA rates.

The models that do not conform to the current CDM-CTOP implementation are still valuable in theoretical analysis, and can be used as benchmarks for the compatible models and as references for future software development.

\section{FCA-PCA Models and PCA Models for FCA Rate Planning}\label{FCA_PCA_Model_section}
In this section, we first focus on two class of models that is proposed for the first subproblem of CTOP traffic rate planning: optimize the acceptance rates given demand estimation. We will use  Enhanced Stochastic Optimization Model (ESOM), which was considered to be the most compatible class of models to the current CDM-CTOP software, as an example to reveal a deficiency that is unique to FCA-PCA models. We also propose a new multistage semi-dynamic FCA-PCA model and argue that a more flexible FCA-PCA model could amplify this shortcoming. We show that how this problem can be addressed in PCA models. The conclusion is that to manage multiple constrained resources precisely and optimally, we need one FCA for one each commodity of traffic flows.

Second, if we further consider flight reroute, because stochastic models are known to give conservative traffic rates, the demand shift issue makes CTOP FCA rate planning very challenging in the stochastic programming framework. This motivates us to also consider to tackle FCA optimization in a different framework, e.g. simulation-based optimization.

\subsection{Enhanced Stochastic Optimization Model (ESOM)}\label{ESOM_discussion}
ESOM extends the classical single resource GDP model to the multiple constrained resources case \cite{HoffmanESOM}. The main advantage of ESOM is that FCAs rates computed from ESOM can be directly implemented in CDM-CTOP software.

For ease of reference, ESOM is listed below. The reader is referred to \cite{HoffmanESOM} for the detailed derivation of this model. Its objective function (\ref{ESOM-Obj}) minimizes ground delay and expected air delay costs; constraints (\ref{ESOM-ground-delay_part1}) and (\ref{ESOM-ground-delay_part2}) plan the ground delay at each FCA; constraint (\ref{ESOM-air-delay}) describes the number of air delayed flights at each PCA; constraint (\ref{ESOM-capacity}) enforces the physical capacity constraint at each PCA; constraint (\ref{ESOM-flow-split}) models the traffic flow split between resources; constraint (\ref{Boundary}) is the boundary condition; constraints (\ref{ESOM-FCA-travel-time}) and (\ref{ESOM-PCA-travel-time}) model the travel time between resources.
\begingroup
\allowdisplaybreaks
\begin{align}
\text{min} \quad & \sum_{r \in \mathcal{F}}\sum_{t\in \mathcal{T}}c_gG_t^r + \sum_{r \in \mathcal{P}}\sum_{t\in \mathcal{T}}\sum_{q\in \mathcal{Q}}c_ap_qA_{t,q}^r \label{ESOM-Obj} \\
 \text{s.t.}\quad  & \tilde{P}^r_t = D^r_t - (G^r_t -G^r_{t-1})& \forall r\in \mathcal{F}, t\in \mathcal{T} \label{ESOM-ground-delay_part1}\\
 &P^r_t = \text{UpFCA}^r_t+ \tilde{P}^r_t & \forall r\in \mathcal{F}, t\in \mathcal{T} \label{ESOM-ground-delay_part2}\\
 & L^r_{t,q} = \text{UpFCA}^r_t + \text{UpPCA}^r_{t,q} -(A^r_{t,q}-A^r_{t-1,q})& \forall r\in \mathcal{P}, q\in\mathcal{Q}, t\in \mathcal{T} \label{ESOM-air-delay}\\
 & M^r_{t,q}\ge L^r_{t,q} & \forall r\in \mathcal{P},  q\in \mathcal{Q}, t\in \mathcal{T} \label{ESOM-capacity}\\
 & \sum_{(r,r^\prime)\in \mathcal{C}}f^{r,r^\prime}_t=1 &\forall r\in \mathcal{R}, t\in \mathcal{T} \label{ESOM-flow-split} \\
 & G^r_t, \tilde{P}^r_t, P^r_t, L^r_{t,q}, A^r_{t,q}\ge 0 & \forall r\in \mathcal{R}, q\in \mathcal{Q}, t\in \mathcal{T} \label{ESOM-Nonnegative} \\
 & G^r_0=G^r_{|\mathcal{T}|}=A^r_{0,q}=A^r_{|\mathcal{T}|,q}=0 & \forall r\in \mathcal{R}, q\in \mathcal{Q}\label{Boundary}\\
 & \text{UpFCA}^r_t = \sum_{(r^\prime, r) \in \mathcal{C}} f^{r^\prime, r}_{t-\Delta^{r^\prime,r}}\cdot P^{r^\prime}_{t-\Delta^{r^\prime,r}} &\forall r\in \mathcal{F}, t\in \mathcal{T} \label{ESOM-FCA-travel-time}\\
 & \text{UpPCA}^r_{t,q} = \sum_{(r^\prime, r) \in \mathcal{C}} f^{r^\prime, r}_{t-\Delta^{r^\prime,r}}\cdot L^{r^\prime}_{t-\Delta^{r^\prime,r},q} &\forall r\in \mathcal{P}, q \in \mathcal{Q}, t\in \mathcal{T} \label{ESOM-PCA-travel-time}
\end{align}
\endgroup
A key trick in ESOM is that \emph{pre-calculated} flow split ratios $f^{r,r^\prime}_t$ are used to model the traffic flow coming out from one resource to multiple downstream resources. We claim as a result, in general we cannot obtain the theoretical optimal FCA planned acceptance rates:
\begin{enumerate}
  \item In ESOM, we cannot follow flights' route schedule precisely. For example, if 10 flights are scheduled to pass PCA1 and then land at the airport in Figure \ref{FCAsPCAs}, in ESOM we cannot guarantee these 10 flights will eventually travel the scheduled planned route and land by the end of planning horizon.

     This is because split ratios $f^{r^\prime, r}_{t}$ are pre-calculated parameters. If we delay some flights to the next time period, we have to use the split ratios in the next period. As an extreme example, in Figure \ref{FCAsPCAs} suppose
    \begin{table}[h!]
      \centering
      \resizebox{0.4\textwidth}{!}{%
      \begin{tabular}{cccccc}
        \toprule
        Time period      & FCA1-PCA1 & FCA1-PCA2 \\
        \midrule
        $t$   & 1        & 0             \\
        $t+1$ & 0        & 1              \\
      \bottomrule
      \end{tabular}}
      \caption{Example of Pre-calculated Ratios}
    \end{table}
     at time period $t$ the split ratio from FCA 1 to PCA1 is 1, and from FCA 1 to PCA 2 is 0; at time period $t+1$, split ratio from FCA 1 to PCA1 is 0, and from FCA 1 to PCA 2 is 1. If we delay some flights from $t$ to $t+1$, then these flights will be rerouted to PCA 2 instead of PCA 1 by ESOM. The flow approximation results in the violation of the original flight route schedules. In some examples, a flight could even be rerouted to a different destination airport.

  \item Since $f^{r^\prime, r}_{t}$ are continuous numbers, in ESOM we cannot impose integrality constraints for decision variables. It is well known that objective value of a linear programming relaxation can be quite different from the objective value of integer solutions \cite{wei2013total}. After we solve ESOM we need to round the planned acceptance rates $P^r_t$ to integers. It is not yet clear what a good rounding strategy is.

  \item In ESOM, we cannot strictly enforce boundary conditions. In TMI rates optimization, we usually require all the affected flights to land at airport or exit the constrained airspace by the end of the planning horizon. This condition is a little tricker to impose in CTOP than in GDP, because in CTOP there can be en route flights between constrained resources. To make sure no flight is still en route at the end, we can choose a sufficiently long planning horizon. Since ground delay and air delay are penalized in the objective function, all CTOP captured flights will eventually exit the FCA-PCA network. However, it is not obvious how long the planning horizon should be. \cite{ZhuAggregate} explicitly enforced that at each PCA the total number of \emph{actual landed/exited flights} during the entire planning horizon should equal to the total number of flights that are \emph{scheduled to land/exit} at that PCA. This approach does not work for ESOM because flow splits are only approximations, the conservation of flow constraints will very likely not be exactly met at each PCA and the problem can become infeasible.
\end{enumerate}
Since traffic flow schedules cannot be strictly followed, boundary conditions may not be exactly met and integrality constraints are not enforced, the objective function is only an approximation to the true objective value.

The root cause of items 1 to 3 is that we approximate the essentially multi-commodity flow as a single-commodity traffic flow. This approximation has to be made if one wants to optimize FCA rates while treating flights as homogeneous traffics.  But in this way, we could only do a compromised traffic flow control and optimization.

\subsection{Multistage Semi-Dynamic ESOM Model}\label{SemiDynamicESOM}
In this section, we introduce the multistage semi-dynamic version of FCA-PCA model, shown in the shaded cell in Table \ref{CTOPTable}. The formulation is listed below:
\begin{align}
\min\quad &\sum_{q\in\mathcal{Q}}p_q \Big(  \sum_{r \in \mathcal{F}} \sum_{s \in \mathcal{S}}\sum_{ \substack{t^\prime=t} }^{|\mathcal{T}|}(t^\prime-t)c_gX^{r,q}_{s,t,t^\prime} + \sum_{r \in \mathcal{P}}\sum_{t\in \mathcal{T}}c_aA^{r}_{t,q}  \Big) \label{Obj_SDESOM}
\end{align}
\begingroup
\allowdisplaybreaks
\begin{align}
&\sum_{ \substack{t^\prime = t}}^{|\mathcal{T}|}X^{r,q}_{s,t,t^\prime} = S^r_{s,t}& \forall s\in \mathcal{S}, t\ge t_s, q\in \mathcal{Q}, r \in \mathcal{F} \label{SD_conservation_of_flow}\\
& P^{r}_{t,q} = \text{UpFCA}^{r}_{t,q} + \sum_{ \substack{ s\in \mathcal{S}\\ t\ge t_s } }\sum_{t\ge t^\prime}X^{r,q}_{s,t^\prime,t} &\forall r \in \mathcal{F}, q\in \mathcal{Q}, t\in \mathcal{T}\label{SD-ESOM-ground-delay} \\
& L^{r}_{t,q} = \text{UpFCA}^{r}_{t,q}+\text{UpPCA}^{r}_{t,q}-(A^{r}_{t,q}-A^{r}_{t-1,q}) & \forall r\in \mathcal{P}, q\in \mathcal{Q}, t\in \mathcal{T} \label{SD_PCA_conservation_of_flow}\\
& L^{r}_{t,q}\le M^{r}_{t,q} & \forall r\in \mathcal{P}, q\in \mathcal{Q}, t\in \mathcal{T}\\
& \sum_{(r,r^\prime)\in \mathcal{C}}f^{r,r^\prime}_t=1 &\quad \forall r\in \mathcal{R}, t\in \mathcal{T} \label{SD_flow_split_ratio}\\
& \text{UpFCA}^{r}_{t,q} = \sum_{(r^\prime, r) \in \mathcal{C}} f^{r^\prime, r}_{t-\Delta^{r^\prime,r}}\cdot P^{r^\prime}_{t-\Delta^{r^\prime,r},q} &\forall r \in \mathcal{F}, q\in \mathcal{Q}, t\in \mathcal{T} \label{SD_UpFCA}\\
& \text{UpPCA}^{r}_{t,q} = \sum_{(r^\prime, r) \in \mathcal{C}} f^{r^\prime, r}_{t-\Delta^{r^\prime,r}}\cdot L^{r^\prime}_{t-\Delta^{r^\prime,r},q} &\forall r \in \mathcal{P}, q\in \mathcal{Q}, t\in \mathcal{T}\\
& A^r_{0,q}=A^r_{|\mathcal{T}|,q}=0 & \forall r\in \mathcal{R}, q\in \mathcal{Q}\\
& X^{r,q}_{s,t,t^\prime}, P^{r}_{t,q}, L^{r}_{t,q}, A^{r}_{t,q}\ge 0& \quad \forall r\in \mathcal{R}, q\in \mathcal{Q}, t\in \mathcal{T}  \label{SD_Nonnegativity}\\
& X^{r,q^b_1}_{s,t,t^\prime} = \dots = X^{r,q^b_{N_b}}_{s,t,t^\prime} &\forall s\in \mathcal{S},t\ge t_s,t^\prime\ge t, b\in \mathcal{B}, t_s\in b \label{SD_Nonanticipativity}
\end{align}
\endgroup
The first set of constraints (\ref{SD_conservation_of_flow}) is the conservation of flow constraints.  At each FCA, we have the planned acceptance rates equal to the summation of the upstream demand and direct demand from the airport (\ref{SD-ESOM-ground-delay}). Constraint (\ref{SD-ESOM-ground-delay}) is similar to (\ref{ESOM-ground-delay_part1})(\ref{ESOM-ground-delay_part2}) in that we only ground-delay flights that directly fly from the departing airports. The difference here is that planned acceptance rates $P^{r}_{t,q}$ now become scenario dependent. As a result, $\text{UpFCA}^r_{t,q}$  is also scenario dependent. Constraints (\ref{SD_PCA_conservation_of_flow})-(\ref{SD_Nonnegativity}) are very similar to (\ref{ESOM-air-delay})-(\ref{ESOM-PCA-travel-time}) in ESOM. In multistage stochastic programming, we need nonanticipativity constraints (\ref{SD_Nonanticipativity}),  which ensure that decisions made at a certain time period are solely based on the information available at that time. Finally, the objective function minimizes expected ground delay and air delay costs.

We want to stress that the semi-dynamic FCA-PCA model has all the imperfections ESOM has. Because of the additional flexibility in this semi-dynamic model, ground delay policy will be differ for different scenarios. Thus, across a range of scenarios, the flow of the air traffic will be even more different than in ESOM. But for all scenarios, $P^{r}_{t,q}$ and $L^{r}_{t,q}$ both use the same split ratio $f^{r,r^\prime}_t$. It is easy to imagine the split ratios may not be good approximations for at least some scenarios. This will be illustrated in a concrete example in section \ref{Model_Comparison}.
\subsection{PCA Models, More FCAs, Optimality}\label{PCA_Optimality}
The primary decision variable in two-stage PCA model is $P^r_{t,\rho}$, which is number of flights belong to path $\rho$ that are planned to be accepted to PCA $r$. \emph{The concept of FCA in TMI corresponds to ground delay assignment}. Thus, we can place one FCA for each path. And FCA designed for controlling flights in one path needs to exempt flights belonging to other paths. For example, in Figure \ref{FCAsPCAs}, there should be actually two FCAs (FCA11 and FCA12) in the position of FCA1, with one controlling flights going through PCA1 and the second one managing flights going through PCA2. Suppose FCA11 is 2 time periods away from PCA1 and $\text{PCA1}\rightarrow \text{PCA3}$ is numbered as path 1, then the acceptance rate of FCA11 at time period $t$ is simply $P^{r=\text{PCA1}}_{t-2,\rho=1}$.

If one decides to  place only one FCA in the location of FCA1, no matter its FCA rates are obtained from ESOM or by adding the rates of FCA11 and FCA12, they are in general NOT optimal.
\subsection{Why Two Subproblems, Demand Shift Issue, Conservative FCA Rates}
The main reason to split CTOP rate planning into two steps is that CTOP slot allocation algorithm is highly nonlinear (section \ref{slot_allocation_algorithm}). Optimizing reroute and delay for each flight while ignoring slot allocation rules is relative easy, as has been done in \cite{ZhuCentralized}. However, it is intractable to optimize reroute and delay at the same time while exactly satisfying all the rules in the slot allocation algorithm.

Demand shift issue appears in the second step of CTOP rate planning problem. It refers to the phenomenon that after we solve the first step problem given a demand estimation and run CTOP slot allocation algorithm, the demand may shift from one FCA to another FCA, and invalidate the proposed acceptance rates. Currently, it is addressed by resolving first step problem whenever demand changes.

Conservative FCA rates issue refers to the fact that in stochastic programming based TMI optimization models, the total number of slots created is equal to total demand. This can be most easily seen in the single airport ground delay model \cite{ball2003stochastic}. If the demand to the airport is zero, then acceptance rates being all zeros will be the optimal solution, even though the airport capacity in the worse scenario can be all greater than zero. If this is considered a minor defect for GDP problem, then it is a major problem in CTOP rate planning. In Figure \ref{FCAsPCAs}, if all flights want to cross PCA2 and land at PCA3, then the planned acceptance rates for PCA1 will be all 0, even though there are capacities in PCA2. CTOP slot allocation algorithm takes these rates as input and will create zero slot for flights to pass through PCA2. Therefore, all flights will continue queuing to pass congested PCA1, which can be very inefficient.

It can be seen that these two issues are coupled: if there is no demand shift issue, even though a conservative FCA rates might not be robust to small demand perturbation, they are acceptable same as in the GDP problem; If FCA rates are not conservative, ideally if they are optimal with respect to any demand information, then demand shift will not be a problem. In \cite{ZhuPARTMI} we have given a counterexample to show that even in the single constrained resource GDP case, in general there does not exist FCA rates that are optimal to any demand. Hence, heuristic need to be used and only suboptimal can be hoped to obtain.

To solve the demand shift issue, we need to have an estimate that how many flights should be accepted to a region of airspace, based only on airspace capacity information. This is not easy because there might be several PCAs connected in different ways and each with very different capacity range. We settle down with finding
approximate upper bounds of FCA acceptance rates, which could achieve the goal that flights can reroute to previously not fully utilized resources. We hope through several iterations of computation, a good planned acceptance rates that enable flights to take better use of all resources can be obtained.


The heuristic we use it called \emph{saturation technique}. The basic idea is that to get approximate upper bounds of FCA rates, we will saturate the FCA-PCA system with artificially high volume of demand. There can be many ways to do it. In this paper, we will test two versions of it:
\begin{enumerate}
  \item We will not use any TOS information and simply apply high demand to all FCAs in each time period. No CTOP slot allocation needs to be run and the computed FCA rates will be the rates we use.
  \item Iterative demand information guided saturation and rate computation
      \begin{enumerate}
    	\item Increase demands to FCAs \emph{proportionally} to sufficiently large numbers  \label{iterative_saturation_step_1}
        \item If the demand to a FCA is 0 at a time period, we will set it as some default small value like 1, such that the planned acceptance rate is not always zero
        \item Run CTOP slot allocation algorithm based on the rates obtained, obtain a new demand estimation, go to step \ref{iterative_saturation_step_1}. Exit if the rates are satisfactory
    \end{enumerate}
\end{enumerate}
The numerical results are shown in section \ref{RCL_comparision}.
\section{Simulation-based Optimization}\label{Simulation_based_Optimization}
Simulation-based optimization integrates computer simulation with optimization technique. It is often applied to problems in which evaluating a solution involves running simulation models. A classical example is aircraft wing design. One has to run complicated computational fluid dynamic models in order to assess a set of proposed parameters \cite{koziel2016simulation}. Simulation-based approach has also been used to in air traffic flow management, for example, in determining the GDP parameters under uncertainly \cite{cook2010model} and in strategically selecting TMI combinations \cite{taylor2015designing}.

The main reason we consider simulation-based optimization is that we can directly optimize FCA rates without worrying about either conservative FCA rates issue or demand shift issue. On the other hand, the CTOP rate planning problem is special compared with \cite{cook2010model} or \cite{taylor2015designing}: in CTOP, we can have up to 5 FCAs, so the dimension of parameter space is much larger than in GDP; to provide real-time decision support, a solution needs to be reported in around 5 minutes, which is much less than in problem \cite{taylor2015designing}. Therefore, tradeoffs and compromises need to be made. In this paper, we use simulation-based optimization as a way to refine the solution we find using other methods. Specifically, we adopt the following two-phase approach:
\begin{enumerate}
  \item In phase one, we will use heuristics including saturation technique to quickly find good initial starting points
  \item In phase two, for a subset of FCAs, we will employ local search method like pattern search to carefully find better solutions
\end{enumerate}
Pattern search is classical derivative free local search method. It is composed of two types of moves: \emph{exploratory search}, which is used to find an improving direction by checking points in the neighbourhood; \emph{pattern move}, which searches in the improving direction and will keep move as long as improvement continues \cite{chinneck2006practical}. In CTOP rate planning problem, all decision variables are integers and bounded, therefore small adjustment about step size and variable range need to be made in the original algorithm.

This two phase framework and pattern search is just one of many options one can choose from. In this paper, we want to use it to obtain a fully CDM-CTOP compatible solution and compare it with various benchmarks in \cite{ZhuAggregate}\cite{ZhuCentralized}. Depending on the complexity of problem in hand and the amount of computing power available, other heuristic algorithms like genetic algorithm and more sophisticatedly design can be considered.

\section{Experimental Results}\label{Experiment}
\begin{figure}
\centering
\begin{minipage}{.5\textwidth}
  \centering
  \includegraphics[width=0.9\linewidth]{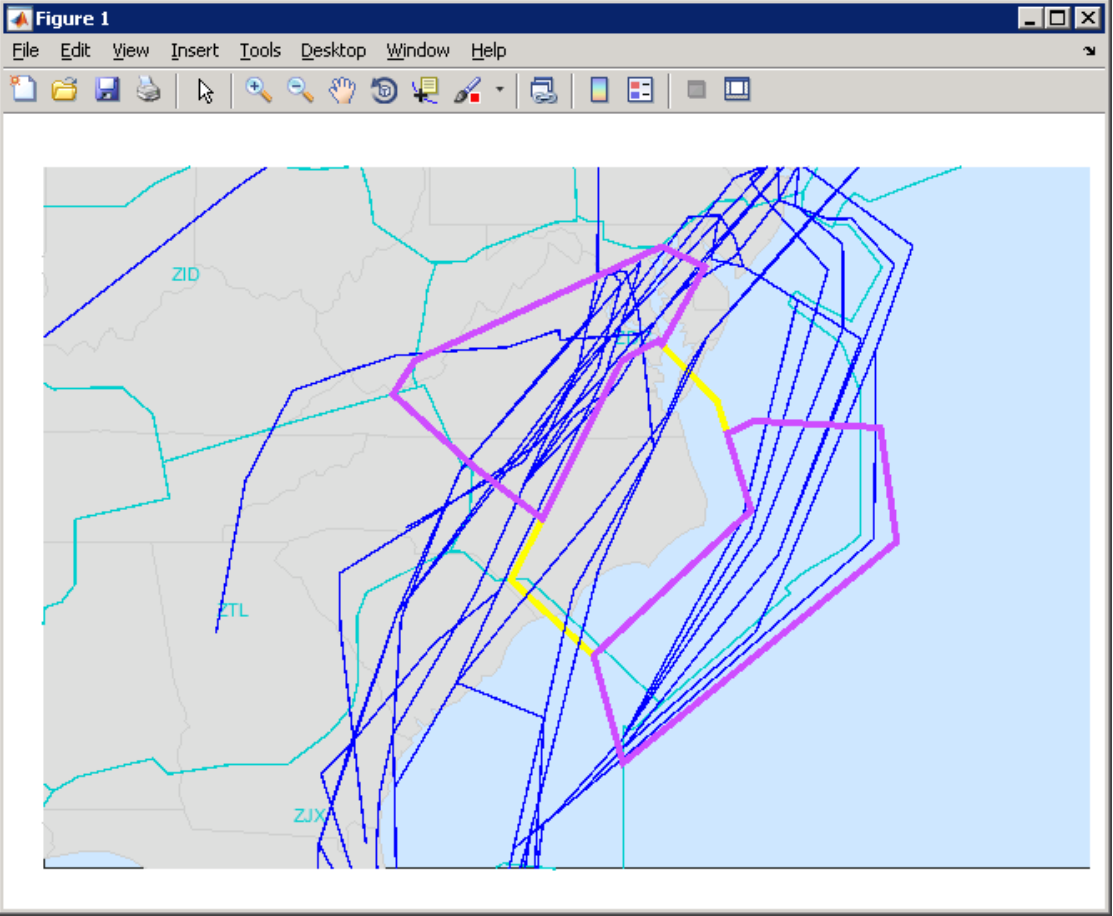}
  \caption{Traffic Trajectory Visualization}\label{Traffic_Routing}
\end{minipage}%
\begin{minipage}{.5\textwidth}
  \centering
  \includegraphics[width=1.03\linewidth]{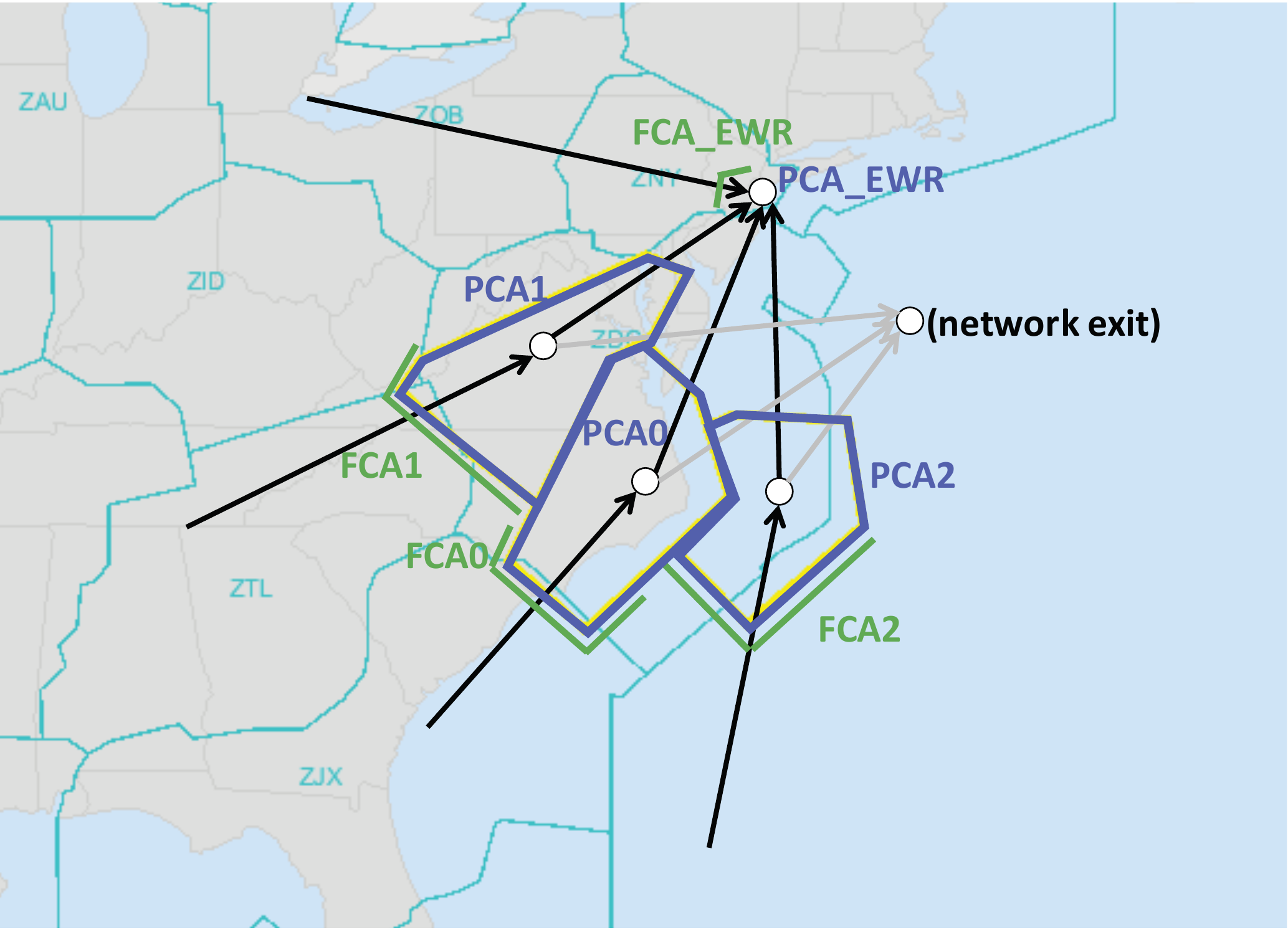}
  \caption{Geographical display of a FCA-PCA network}\label{PCA_Network}
\end{minipage}
\end{figure}

To compare the performance of the proposed semi-dynamic model with ESOM and PCA models, we continue using the test case with convective weather activity in southern Washington Center (ZDC) and EWR airport.  We assume there is a four-hour capacity reduction in ZDC/EWR from 2000Z to 2359Z. By analyzing the traffic trajectory and weather data, we can build the FCA-PCA network, shown in Figure \ref{PCA_Network}. In this use case, each FCA directly lies atop of the corresponding PCA. As we discussed section \ref{PCA_Optimality}, for example ideally there should be two FCAs in front of PCA1, instead of just FCA1. However, the FCA-PCA network in Figure \ref{PCA_Network} is what most air traffic manager would design and it is what our subject matter experts recommended. We will see how various models perform on this network.

All optimization models are solved using Gurobi 8.1 on a workstation with 3.6 GHz processors and 32 GB RAM.


\subsection{Capacity Profile and Traffic Demand}
In this work, we directly manipulate the capacity profiles from the base forecast to create alternative capacity profiles. This gives us full control over the capacity data for experimental purposes. The capacity information is shown in Table \ref{Capacity_Scenarios}. We can see that in scenario 1 at 2100Z PCA1's 15-minute capacity changes from 44 to 50, the EWR's capacity changes from 8 to 10; in scenario 2 at 2230Z, the capacities of PCA1 and EWR return to the nominal values.  These two changes correspond to the two branch points in the scenario tree shown in Figure \ref{ScenarioTreeExperiment}.

In GDP optimization, we usually add one extra time period to make sure all flights will land at the end of the planning horizon. Because CTOP has multiple constrained resources, we need to add more than one time period depending on the topology of the FCA-PCA network. In this case, we add eight extra time periods to make sure all flights will land/exit by the end of planning horizon. For any time periods outside the CTOP start-end time, e.g. the extra eight time periods in Table \ref{Capacity_Scenarios}, we assumed nominal capacity.

\begin{table*}[t]
  \resizebox{1.0\textwidth}{!}{%
  \centering
  \begin{tabular}{cc >{\columncolor[gray]{0.8}}c >{\columncolor[gray]{0.8}}c >{\columncolor[gray]{0.8}} c >{\columncolor[gray]{0.8}}c cccc
  >{\columncolor[gray]{0.8}}c >{\columncolor[gray]{0.8}}c >{\columncolor[gray]{0.8}}c >{\columncolor[gray]{0.8}}c
  cccc
  >{\columncolor[gray]{0.8}}c  >{\columncolor[gray]{0.8}}c >{\columncolor[gray]{0.8}} c >{\columncolor[gray]{0.8}}c c ccc }
  \toprule
  & Resource & 20:00 & 15 & 30 & 45 & 21:00 & 15 & 30 & 45 & 22:00 & 15& 30 & 45 & 23:00 & 15 & 30 & 45 & 00:00 & 15 & 30 & 45 &01:00&15&30&45\\
  \midrule
  \multirow{4}{*}{Scen1} & PCA0  & 13 & 13 & 13 & 13 & 25 & 25 & 25 & 25 & 25 & 25 & 25 & 25 & 25 & 25 & 25 & 25 & 25 & 25 & 25 & 25&25& 25 & 25 & 25\\
                         & PCA1 & 44 & 44 & 44 & 44 & 50 & 50 & 50 & 50 & 50 & 50 & 50 & 50 & 50 & 50 & 50 & 50 & 50 & 50 & 50 & 50 &50& 50 & 50 & 50\\
                         & PCA2 &  5 &  5 &  5 &  5 &  5 &  5 &  5 &  5 &  5 &  5 &  5 &  5 &  5 &  5 &  5 &  5 &  5 &  5 &  5 &  5 & 5&  5 &  5 &  5\\
                         & EWR  &  8 &  8 &  8 &  8 & 10 & 10 & 10 & 10 & 10 & 10 & 10 & 10 & 10 & 10 & 10 & 10 & 10 & 10 & 10 & 10  & 10 & 10 & 10 & 10\\
  \midrule
  \multirow{4}{*}{Scen2} & PCA0  & 13 & 13 & 13 & 13 & 13 & 13 & 13 & 13 & 13 & 13 & 25 & 25 & 25 & 25 & 25 & 25	& 25 & 25 & 25 & 25 & 25& 25 & 25 & 25\\
                         & PCA1 & 44 & 44 & 44 & 44 & 44 & 44 & 44 & 44 & 44 & 44 & 50 & 50 & 50 & 50 & 50 & 50 & 50 & 50 & 50 & 50 & 50 & 50 & 50 & 50\\
                         & PCA2 &  5 &  5 &  5 &  5 &  5 &  5 &  5 &  5 &  5 &  5 &  5 &  5 &  5 &  5 &  5 &  5 &  5 &  5 &  5 &  5 & 5&  5 &  5 &  5 \\
                         & EWR  &  8 &  8 &  8 &  8 &  8 &  8 &  8 &  8 &  8 &  8 & 10 & 10 & 10 & 10 & 10 & 10 & 10 & 10 & 10 & 10 & 10& 10 & 10 & 10 \\
  \midrule
  \multirow{4}{*}{Scen3} & PCA0  & 13 & 13 & 13 & 13 & 13 & 13 & 13 & 13 & 13 & 13 & 13 & 13 & 13 & 13 & 13 & 13 & 25 & 25 & 25 & 25 & 25& 25 & 25 & 25\\
                         & PCA1 & 44 & 44 & 44 & 44 & 44 & 44 & 44 & 44 & 44 & 44 & 44 & 44 & 44 & 44 & 44 & 44 & 50 & 50 & 50 & 50 & 50& 50 & 50 & 50\\
                         & PCA2 &  5 &  5 &  5 &  5 &  5 &  5 &  5 &  5 &  5 &  5 &  5 &  5 &  5 &  5 &  5 &  5 &  5 &  5 &  5 &  5 & 5&  5 &  5 &  5\\
                         & EWR  &  8 &  8 &  8 &  8 &  8 &  8 &  8 &  8 &  8 &  8 &  8 &  8 &  8 &  8 &  8 &  8 & 10 & 10 & 10 & 10 & 10& 10 & 10 & 10 \\
  \bottomrule
  \end{tabular}}
  \caption{Capacity Scenarios} \label{Capacity_Scenarios}
\end{table*}

We use historical flight data for traffic demand modeling. We only keep flights which pass through one of the 3 PCAs created in ZDC plus all EWR arrivals. The resulting set contains 1098 flights; among them, 890 flights traverse the PCAs in their active periods and 130 flights land at EWR airport. There are 1368 TOS options for these 890 flights, on average 1.54 options per flight.


\subsection{Stochastic Model Comparisons in terms of System Delay Costs}\label{Model_Comparison}
\begin{table*}[t]
  \resizebox{1.0\textwidth}{!}{%
  \centering
    \begin{tabular}{c c c c c c c c c c c c}
    \toprule
     \multirow{3}{*}{} & \multicolumn{3}{c}{Ground Delay Periods} & \multicolumn{3}{c}{Air Holding Periods} & \multicolumn{3}{c}{Total Cost} & \multirow{3}{*}{ Expected Cost } & \multirow{3}{*}{ Running Time }\\
      & \multicolumn{3}{c}{If This Scenario Occurs:} & \multicolumn{3}{c}{If This Scenario Occurs:} & \multicolumn{3}{c}{If This Scenario Occurs:} & \\
      \cmidrule{2-10}
     &SCEN1 & SCEN2 & SCEN3  &SCEN1 & SCEN2 & SCEN3        &SCEN1 & SCEN2 & SCEN3  & &Seconds\\
    \midrule
     ESOM                                  & 296.05 & 296.05 & 296.05 & 0 & 0 & 211.55 & 296.05 & 296.05 & 719.15 & {\cellcolor[gray]{0.8}} 422.98 & 0.03\\
     Semi-Dynamic ESOM                     & 93.80  & 292.13 & 507.60 & 0 & 0 & 0      & 93.80  & 292.13 & 507.60 & {\cellcolor[gray]{0.8}} 297.27 & 0.22\\
     \midrule
     Two-Stage Model \cite{ZhuAggregate}   & 284 & 284 & 284  & 0 & 0  & 200 & 284 & 284 & 684  & {\cellcolor[gray]{0.8}}404.0& 0.17\\
     Semi-Dynamic Model \cite{ZhuAggregate}& 165 & 284 & 415  & 0 & 0  &  69 & 165 & 284 & 484  & {\cellcolor[gray]{0.8}}329.0& 0.67\\
     Dynamic Model \cite{ZhuAggregate}     & 126 & 284 & 479  & 0 & 0  &  9  & 126 & 284 & 488  & {\cellcolor[gray]{0.8}}300.5& 1.19\\
  \bottomrule
  \end{tabular}}
  \caption{ FCA-PCA vs. PCA Models Stochastic Solutions Comparison}\label{TwoStageTable}
\end{table*}

\begin{table*}[t]
  \resizebox{1.0\textwidth}{!}{%
  \centering
  \begin{tabular}{ccc c  c c cccc
  c c c c
  cccc
  c  c  c c c ccc }
  \toprule
  & Resource & 20:00 & 15 & 30 & 45 & 21:00 & 15 & 30 & 45 & 22:00 & 15& 30 & 45 & 23:00 & 15 & 30 & 45 & 00:00 & 15 & 30 & 45 &01:00&15&30&45\\
  \toprule
    \multirow{3}{*}{\pbox{3cm}{Flow Split\\Ratio EWR}}
    & FCA0   & 0      & 0.04  & 0.043 & 0     & 0     & 0     & 0.083 & 0.133 & 0 & 0 & 0     & 0     & 0   & 0     & 0     & 0 & 0 & 0 & 0 & 0 & 0.5& 0.5& 0.5& 0.5\\
    & FCA1   & 0.025  & 0.064 & 0.095 & 0.233 & 0.067 & 0.023 & 0.1   & 0.065 & 0.105 & 0.031 & 0.045 & 0.1 & 0.067 & 0.125 & 0.057 & 0.053 & 0.121 & 0.133 & 0.074 & 0.057 & 0.5& 0.5& 0.5& 0.5\\
    & FCA2   & 0.5    & 0.5   &     0 &   0.5 & 0.5   & 0.5   & 0.5   & 0.5   & 0.5   & 0.5   & 0     & 0.5 & 0.5   & 0.5   & 0.5   & 0.5   & 0.5   & 0.5   & 0.5   & 0.5   & 0.5& 0.5& 0.5& 0.5\\
   \toprule
    \multirow{4}{*}{\pbox{3cm}{ESOM\\Model}}        & FCA0     & 13 & 13 & 13 & 13 & 13 & 13 & 13 & 13 & 13 & 13 & 25 & 25 & 22 & 12 & 12 &  8 & 0 & 0 & 0  & 0 & 0 & 0 & 0 & 0\\
                                                    & FCA1     & 40 & 44 & 44 & {\cellcolor[gray]{0.8}}22 & {\cellcolor[gray]{0.8}}39 & 44 & 30 & 44 & 40 & 32 & 44 & 30 & 30 & 32 & 35 & 19 & 0 & 0 & 0  & 0 & 0 & 0 & 0 & 0\\
                                                    & FCA2     &  0 &  0 &  1 &  0 &  0 &  0 &  0 &  0 &  0 &  0 &  1 &  0 &  0 &  0 &  0 &  0 & 0 & 0 & 0 & 0  & 0 & 0 & 0 & 0\\
                                                    & FCA-EWR  &  8 &  8 &  5 &  5 &  4 &  2 &  {\cellcolor[gray]{0.8}}3 &  {\cellcolor[gray]{0.8}}2 &  {\cellcolor[gray]{0.8}}5 &  7 &  5 &  5 &  {\cellcolor[gray]{0.8}}4 &  {\cellcolor[gray]{0.8}}9 &  {\cellcolor[gray]{0.8}}8 &  4 & 0 & 0 & 0 & 0  & 0 & 0 & 0 & 0\\
  \midrule
  \multirow{4}{*}{\pbox{3cm}{Two-stage\\PCA Model}} & PCA0     & 13 & 13 & 13 & 13 & 13 & 13 & 13 & 13 & 13 & 13 & 25 & 25 & 22 & 12 & 12 &  8 & 0 & 0 & 0  & 0 & 0 & 0 & 0 & 0\\
                                                    & PCA1     & 40 & 44 & 44 & {\cellcolor[gray]{0.8}}31 & {\cellcolor[gray]{0.8}}30 & 44 & 30 & 44 & 40 & 32 & 44 & 30 & 30 & 32 & 35 & 19 & 0 & 0 & 0 & 0  & 0 & 0 & 0 & 0\\
                                                    & PCA2     &  0 &  0 &  1 &  0 &  0 &  0 &  0 &  0 &  0 &  0 &  1 &  0 &  0 &  0 &  0 &  0 & 0 & 0 & 0 & 0  & 0 & 0 & 0 & 0\\
                                                    & EWR  &  8 &  8 &  5 &  5 &  4 &  2 &  {\cellcolor[gray]{0.8}}4 &  {\cellcolor[gray]{0.8}}1 &  {\cellcolor[gray]{0.8}}6 &  7 &  5 &  5 &  {\cellcolor[gray]{0.8}}6 &  {\cellcolor[gray]{0.8}}8 &  {\cellcolor[gray]{0.8}}7 &  4 & 0 & 0 & 0 & 0  & 0 & 0 & 0 & 0\\
  \midrule
  \end{tabular}}
  \caption{ FCA-PCA vs. PCA Two-stage Models Acceptance Rate Comparison}\label{ESOM_Rate_Comparision}
\end{table*}

\begin{table*}[t]
  \centering
  \resizebox{0.9\textwidth}{!}{%
  \begin{tabular}{ccc c  c c cccc
  c c c c
  cccc }
  \toprule
  & Resource & 20:00 & 15 & 30 & 45 & 21:00 & 15 & 30 & 45 & 22:00 & 15& 30 & 45 & 23:00 & 15 & 30 & 45 \\
   \toprule
    \multirow{4}{*}{\pbox{3cm}{ESOM\\Model}}        & FCA0     & 0 & 12 & 22 & 25 & 31 & 28 & 27 & 29 & 29 & 30 & 25 & 11 & 0 & 0 & 0 & 0 \\
                                                    & FCA1     & 0 &  3 &  1 &  {\cellcolor[gray]{0.8}}9 & 0  &  0 &  0 & 2  &  0 &  0 &  0 &  0 & 0 & 0 & 0 & 0 \\
                                                    & FCA2     & 0 &  0 &  0 & 0  &  0 &  0 &  0 & 0  &  0 &  0 &  0 &  0 & 0 & 0 & 0 & 0 \\
                                                    & FCA-EWR  & 0 &  0 &  0 & 0  &  0 &  0 &  1 & 3  &  2 &  1 &  0 &  0 & 4 & 1 & 0 & 0 \\
  \midrule
  \multirow{4}{*}{\pbox{3cm}{Two-stage\\PCA Model}} & PCA0     & 0 & 12 & 22 & 25 & 31 & 28 & 27 & 29 & 29 & 30 & 25 & 11 & 0 & 0 & 0 & 0 \\
                                                    & PCA1     & 0 &  3 &  1 & {\cellcolor[gray]{0.8}}0  &  0 &  0 &  0 & 2  &  0 &  0 &  0 &  0 & 0 & 0 & 0 & 0 \\
                                                    & PCA2     & 0 &  0 &  0 & 0  &  0 &  0 &  0 & 0  &  0 &  0 &  0 &  0 & 0 & 0 & 0 & 0 \\
                                                    & EWR      & 0 &  0 &  0 & 0  &  0 &  0 &  0 & 4  &  2 &  1 &  0 &  0 & 2 & 0 & 0 & 0 \\
  \midrule
  \end{tabular}}
  \caption{ FCA-PCA vs. PCA Two-stage Models Ground Delay Comparison}\label{Ground_Delay_Comparision}
\end{table*}

\begin{table*}[t]
  \resizebox{1.0\textwidth}{!}{%
  \centering
  \begin{tabular}{ccccccc cccc
  cccc
  cccc
  cccc cccc }
  \toprule
  Models&Scenario& Resource & 20:00 & 15 & 30 & 45 & 21:00 & 15 & 30 & 45 & 22:00 & 15& 30 & 45 & 23:00 & 15 & 30 & 45 & 00:00 & 15 & 30 & 45 &01:00&15&30&45\\
  \midrule
  \multirow{12}{*}{ \pbox{3cm}{Semi-Dynamic\\ESOM}} &\multirow{4}{*}{Scen1}
                          & PCA0 & 13 & 13 & 13 & 13 & {\cellcolor[gray]{0.8}}25 & {\cellcolor[gray]{0.8}}25 & {\cellcolor[gray]{0.8}}16 & {\cellcolor[gray]{0.8}}15 & {\cellcolor[gray]{0.8}}13 & 14 & {\cellcolor[gray]{0.8}}20 & {\cellcolor[gray]{0.8}}11 & 11 & 12 & 12 & 8  & 0 & 0 & 0 & 0 & 0 & 0 & 0 & 0  \\
                         && PCA1 & 40 & 44 & 44 & {\cellcolor[gray]{0.8}}31 & {\cellcolor[gray]{0.8}}30 & 44 & 30 & 46 & 38 & 32 & 44 & 30 & 30 & 32 & 35 & 19 & 0 & 0 & 0 & 0 & 0 & 0 & 0 & 0  \\
                         && PCA2 &  0 &  0 &  1 &  0 &  0 &  0 &  0 &  0 &  0 &  0 &  1 &  0 &  0 &  0 &  0 &  0 & 0 & 0 & 0 & 0  & 0 & 0 & 0 & 0\\
                         && EWR  & 8 & 8 & 5 & 5 & 4 & 2 & 4 & {\cellcolor[gray]{0.8}}2 & {\cellcolor[gray]{0.8}}7 & {\cellcolor[gray]{0.8}}6 & 4 & 5 & {\cellcolor[gray]{0.8}}4 & {\cellcolor[gray]{0.8}}9 & {\cellcolor[gray]{0.8}}8 & 4 & 0 & 0 & 0 & 0 & 0 & 0 & 0 & 0\\
  \cmidrule[\heavyrulewidth]{2-27}
  &\multirow{4}{*}{Scen2} & PCA0 & 13 & 13 & 13 & 13 & 13 & 13 & 13 & 13 & 13 & 13 & 25 & 25 & 22 & 12 & 12 & 8  & 0 & 0 & 0 & 0 & 0 & 0 & 0 & 0\\
                         && PCA1 & 40 & 44 & 44 & {\cellcolor[gray]{0.8}}31 & {\cellcolor[gray]{0.8}}30 & 44 & 30 & 44 & 40 & 32 & 44 & 30 & 30 & 32 & 35 & 19 & 0 & 0 & 0 & 0 & 0 & 0 & 0 & 0 \\
                         && PCA2 &  0 &  0 &  1 &  0 &  0 &  0 &  0 &  0 &  0 &  0 &  1 &  0 &  0 &  0 &  0 &  0 & 0 & 0 & 0 & 0  & 0 & 0 & 0 & 0\\
                         && EWR  & 8 & 8 & 5 & 5 & 4 & 2 & {\cellcolor[gray]{0.8}}3 & {\cellcolor[gray]{0.8}}0 & {\cellcolor[gray]{0.8}}6 & 7 & {\cellcolor[gray]{0.8}}7 & 5 & {\cellcolor[gray]{0.8}}4 & {\cellcolor[gray]{0.8}}9 & {\cellcolor[gray]{0.8}}8 & {\cellcolor[gray]{0.8}}4 & 0 & 0 & 0 & 0 & 0 & 0 & 0 & 0\\
  \cmidrule[\heavyrulewidth]{2-27}
  &\multirow{4}{*}{Scen3} & PCA0 & 13 & 13 & 13 & 13 & 13 & 13 & 13 & 13 & 13 & 13 & {\cellcolor[gray]{0.8}}13 & {\cellcolor[gray]{0.8}}13 & {\cellcolor[gray]{0.8}}13 & {\cellcolor[gray]{0.8}}13 & {\cellcolor[gray]{0.8}}13 & {\cellcolor[gray]{0.8}}13 & 8 & 6 & 12 & 0 & 0 & 0 & 0 & 0 \\
                         && PCA1 & 40 & 44 & 44 & {\cellcolor[gray]{0.8}}22 & {\cellcolor[gray]{0.8}}39 & 44 & 30 & 44 & 40 & 32 & {\cellcolor[gray]{0.8}}44 & {\cellcolor[gray]{0.8}}30 & {\cellcolor[gray]{0.8}}30 & 32 & 35 & 19 & 0 & 0 & 0  & 0 & 0 & 0 & 0 & 0 \\
                         && PCA2 &  0 &  0 &  1 &  0 &  0 &  0 &  0 &  0 &  0 &  0 &  1 &  0 &  0 &  0 &  0 &  0 & 0 & 0 & 0  & 0  & 0 & 0 & 0 & 0\\
                         && EWR  & 8 & 8 & 5 & 5 & 4 & 2 & {\cellcolor[gray]{0.8}}3 & {\cellcolor[gray]{0.8}}2 & {\cellcolor[gray]{0.8}}5 & 7 & 5 & {\cellcolor[gray]{0.8}}4 & {\cellcolor[gray]{0.8}}2 & 7 & {\cellcolor[gray]{0.8}}6 & {\cellcolor[gray]{0.8}}5 & {\cellcolor[gray]{0.8}}4 & {\cellcolor[gray]{0.8}}2 & 0 & 0 & 0 & 0 & 0 & 0\\
  \midrule
  \midrule
    \multirow{12}{*}{\pbox{3cm}{Semi-Dynamic\\PCA}} &\multirow{4}{*}{Scen1}
                          & FCA0 & 13 & 13 & 13 & 13 & {\cellcolor[gray]{0.8}}15 & {\cellcolor[gray]{0.8}}15 & {\cellcolor[gray]{0.8}}17 & {\cellcolor[gray]{0.8}}22 & {\cellcolor[gray]{0.8}}20 & 14 & {\cellcolor[gray]{0.8}}23 & {\cellcolor[gray]{0.8}}13 & 11 & 12 & 12 & 8  & 0 & 0 &  0 & 0 & 0 & 0 & 0 & 0 \\
                         && PCA1 & 40 & 44 & 44 & {\cellcolor[gray]{0.8}}29 & {\cellcolor[gray]{0.8}}32 & 44 & 30 & 46 & 38 & 32 & 44 & 30 & 30 & 32 & 35 & 19 & 0 & 0 &  0 & 0 & 0 & 0 & 0 & 0 \\
                         && PCA2 &  0 &  0 &  1 &  0 &  0 &  0 &  0 &  0 &  0 &  0 &  1 &  0 &  0 &  0 &  0 &  0 & 0 & 0 & 0 & 0  & 0 & 0 & 0 & 0\\
                         && EWR  & 8 & 8 & 5 & 5 & 4 & 2 & 4 & {\cellcolor[gray]{0.8}}4 & {\cellcolor[gray]{0.8}}4 & {\cellcolor[gray]{0.8}}7 & 4 & 5 & {\cellcolor[gray]{0.8}}6 & {\cellcolor[gray]{0.8}}8 & {\cellcolor[gray]{0.8}}7 & 4 & 0 & 0 & 0 & 0 & 0 & 0 & 0 & 0  \\
  \cmidrule[\heavyrulewidth]{2-27}
  &\multirow{4}{*}{Scen2} & PCA0 & 13 & 13 & 13 & 13 & 13 & 13 & 13 & 13 & 13 & 13 & 25 & 25 & 22 & 12 & 12 & 8  & 0 & 0 &  0 & 0 & 0 & 0 & 0  & 0 \\
                         && PCA1 & 40 & 44 & 44 & {\cellcolor[gray]{0.8}}29 & {\cellcolor[gray]{0.8}}32 & 44 & 30 & 44 & 40 & 32 & 44 & 30 & 30 & 32 & 35 & 19 & 0 & 0 &  0 & 0 & 0 & 0 & 0  & 0 \\
                         && PCA2 &  0 &  0 &  1 &  0 &  0 &  0 &  0 &  0 &  0 &  0 &  1 &  0 &  0 &  0 &  0 &  0 & 0 & 0 & 0 & 0  & 0 & 0 & 0 & 0\\
                         && EWR  & 8 & 8 & 5 & 5 & 4 & 2 & {\cellcolor[gray]{0.8}}4 & {\cellcolor[gray]{0.8}}3 & {\cellcolor[gray]{0.8}}4 & 7 & {\cellcolor[gray]{0.8}}5 & 5 & {\cellcolor[gray]{0.8}}6 & {\cellcolor[gray]{0.8}}8 & {\cellcolor[gray]{0.8}}7 & {\cellcolor[gray]{0.8}}4 & 0 & 0 & 0 & 0 & 0 & 0 & 0 & 0  \\
  \cmidrule[\heavyrulewidth]{2-27}
  &\multirow{4}{*}{Scen3} & PCA0 & 13 & 13 & 13 & 13 & 13 & 13 & 13 & 13 & 13 & 13 & {\cellcolor[gray]{0.8}}24 & {\cellcolor[gray]{0.8}}20 & {\cellcolor[gray]{0.8}}17 &  {\cellcolor[gray]{0.8}}4 &  {\cellcolor[gray]{0.8}}1 & {\cellcolor[gray]{0.8}}12 & 8 & 6 & 12 & 0 & 0 & 0 & 0 & 0 \\
                         && PCA1 & 40 & 44 & 44 & {\cellcolor[gray]{0.8}}29 & {\cellcolor[gray]{0.8}}32 & 44 & 30 & 44 & 40 & 32 &{\cellcolor[gray]{0.8}}43 &{\cellcolor[gray]{0.8}}29 &{\cellcolor[gray]{0.8}}32 & 32 & 35 & 19 & 0 & 0 &  0 & 0 & 0 & 0 & 0 & 0 \\
                         && PCA2 &  0 &  0 &  1 &  0 &  0 &  0 &  0 &  0 &  0 &  0 &  1 &  0 &  0 &  0 &  0 &  0 & 0 & 0 & 0 & 0  & 0 & 0 & 0 & 0\\
                         && EWR  & 8 & 8 & 5 & 5 & 4 & 2 & {\cellcolor[gray]{0.8}}4 & {\cellcolor[gray]{0.8}}3 & {\cellcolor[gray]{0.8}}4 & 7 & 5 & {\cellcolor[gray]{0.8}}5 & {\cellcolor[gray]{0.8}}6 & 7 & {\cellcolor[gray]{0.8}}5 & {\cellcolor[gray]{0.8}}6 & {\cellcolor[gray]{0.8}}1 & {\cellcolor[gray]{0.8}}0 & 0 & 0 & 0 & 0 & 0 & 0 \\
  \bottomrule
  \end{tabular}}
  \caption{FCA-PCA vs. PCA Semi-dynamic Models Acceptance Rate Comparison}\label{SD_Rate_Comparision}
\end{table*}

\begin{figure}[h!]
\centering
\begin{minipage}{.5\textwidth}
  \centering
  \includegraphics[width=\linewidth]{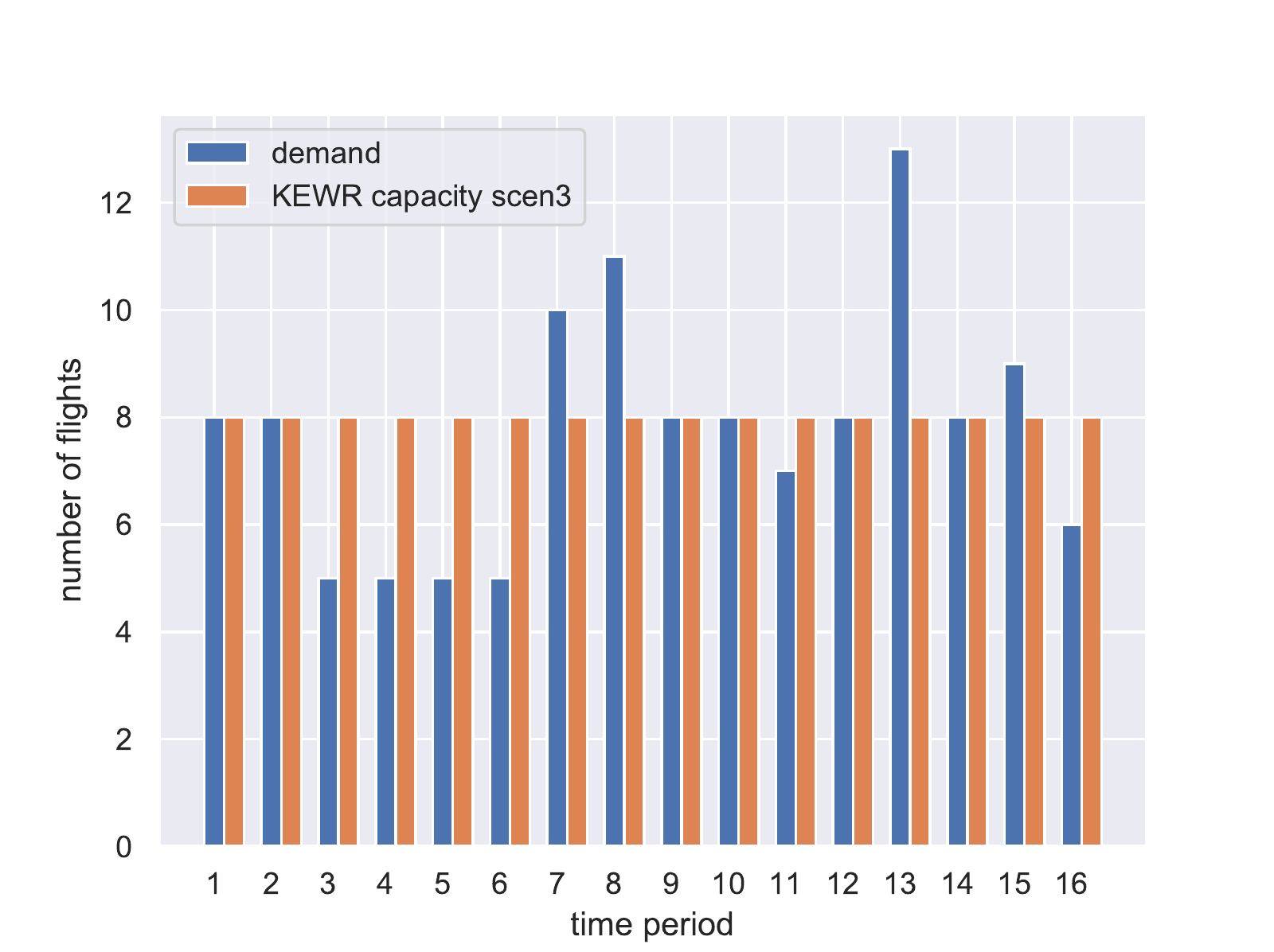}
  \caption{Demand and Capacity at EWR}\label{EWR_demand_capacity_fig}
\end{minipage}%
\begin{minipage}{.5\textwidth}
  \centering
  \includegraphics[width=\linewidth]{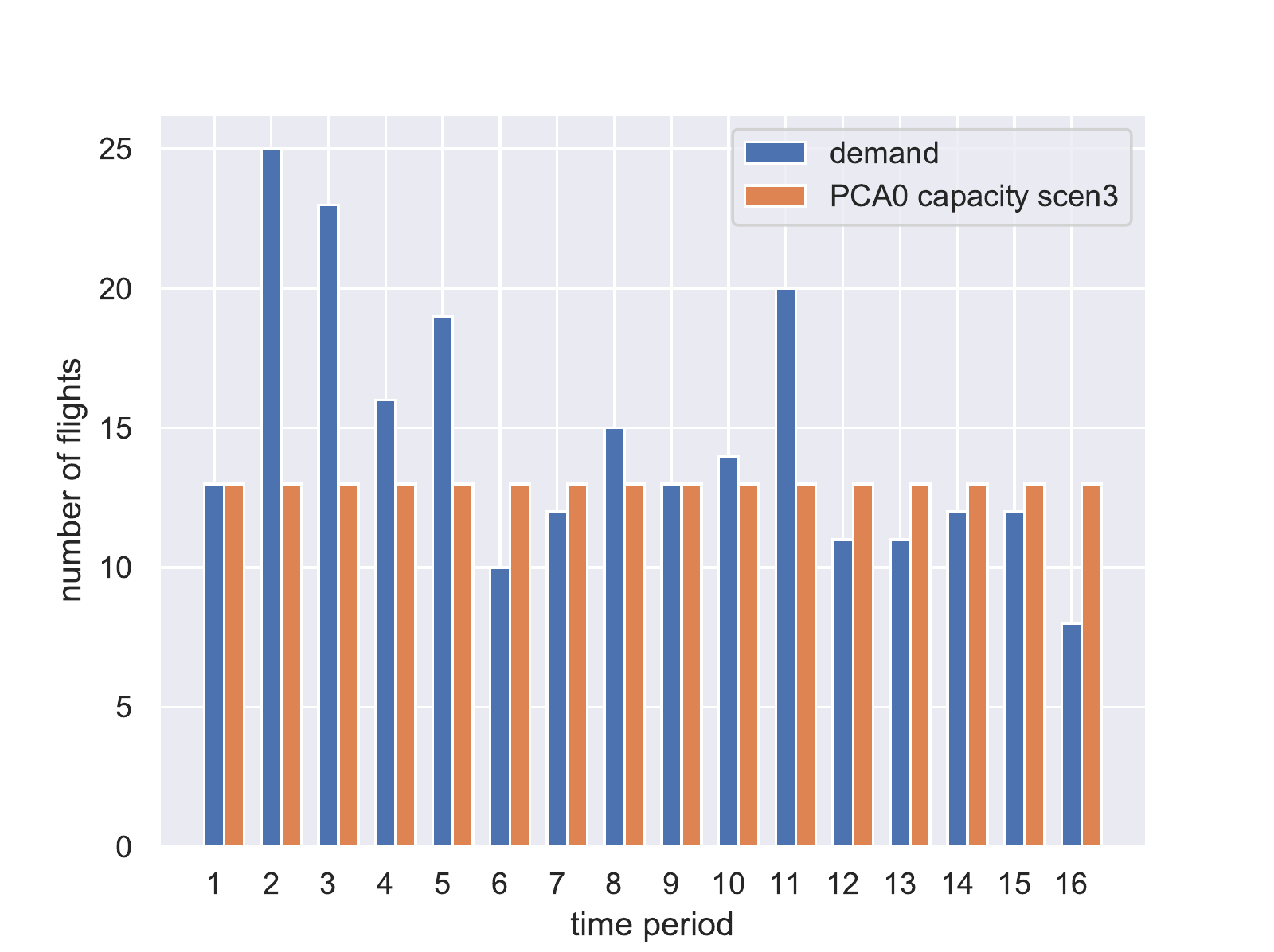}
  \caption{Demand and Capacity at PCA0}
\end{minipage}
\begin{minipage}{.5\textwidth}
  \centering
  \includegraphics[width=\linewidth]{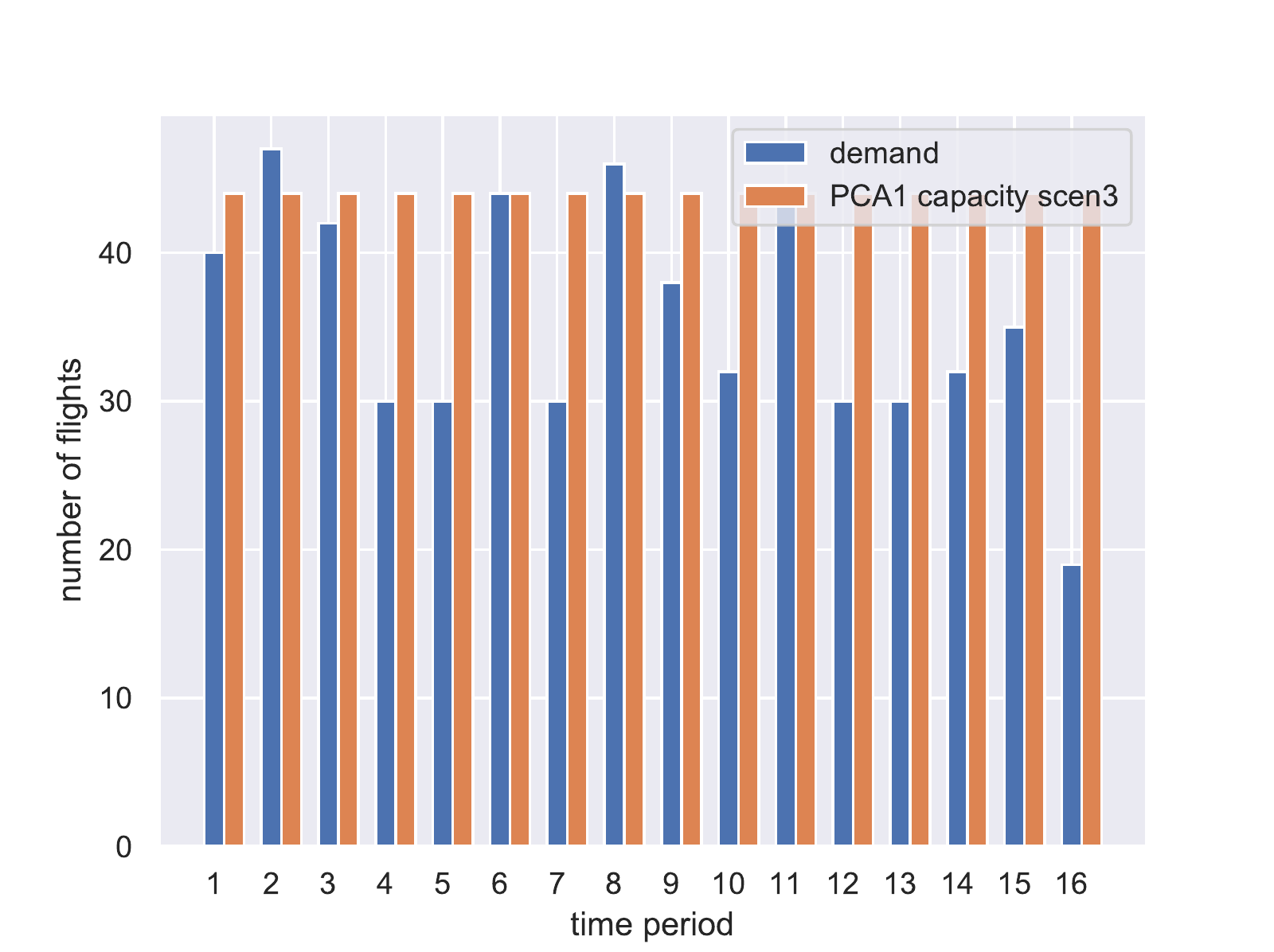}
  \caption{Demand and Capacity at PCA1}
\end{minipage}%
\begin{minipage}{.5\textwidth}
  \centering
  \includegraphics[width=\linewidth]{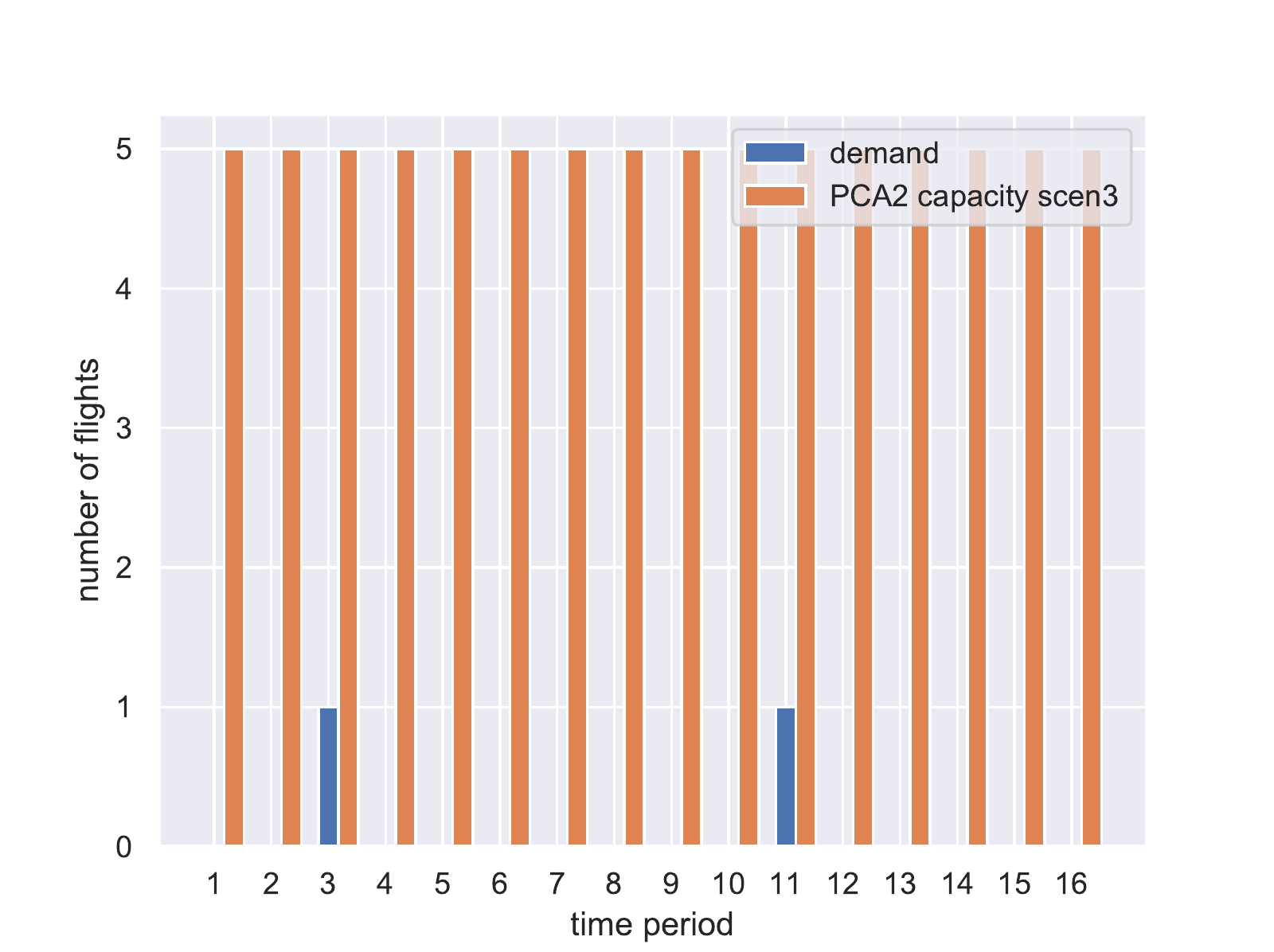}
  \caption{Demand and Capacity at PCA2}
\end{minipage}
\end{figure}
The main results are listed in Table \ref{TwoStageTable}-\ref{SD_Rate_Comparision}. In these four tables, for PCA models, the rates and ground delay at each PCA is the summation of rates and ground delay contributed by all paths that contains that PCA. There are some key observations from this table: 
\begin{enumerate}
  \item Even though we did not enforce integrality constraints on decision variables in ESOM, still the objective value of ESOM is larger than the corresponding two-stage integer PCA model. This unexpected result is exactly due to items 1 to 3 discussed in section \ref{ESOM_discussion}. For example, in Table \ref{Ground_Delay_Comparision}, we can see the ground delay for FCA1 at 20:45 is 9 for ESOM, which is much larger than 0 as determined by two-stage PCA model. This is because flow split ratio for FCA1 to EWR at 20:45 is 0.233, which is much larger than 0.067 at 21:00. EWR is rather congested after 21:30 (time period 7), as can be seen in Figure \ref{EWR_demand_capacity_fig}. By accepting more flights at 21:00, lower flow split ratio is exploited and fewer flights will go to EWR and take more expensive air delay.
  \item The objective value of semi-dynamic ESOM is smaller than semi-dynamic PCA model. Even though the result is the other way around, the underlying explanations are similar. If we delay some flights to a time period in which the fraction of traffic to the constrained resources is small, some of the demands essentially vanish. That is the one of the reasons we see a lot of ground delays in scenario 3 for semi-dynamic ESOM.
  \item In this example, the more flexible the model is, the larger the difference between FCA-PCA model and PCA model in terms of objective value (Table \ref{TwoStageTable}) and planned acceptance rate (Table \ref{ESOM_Rate_Comparision} and Table \ref{SD_Rate_Comparision}). This is because a more flexible model can more heavily exploit the imperfections of the FCA-PCA model.
\end{enumerate}

\subsection{Stochastic and Simulation-based Optimization Result}\label{RCL_comparision}
\begin{figure}[h!]
\centering
\begin{minipage}{.5\textwidth}
  \centering
  \includegraphics[width=\linewidth]{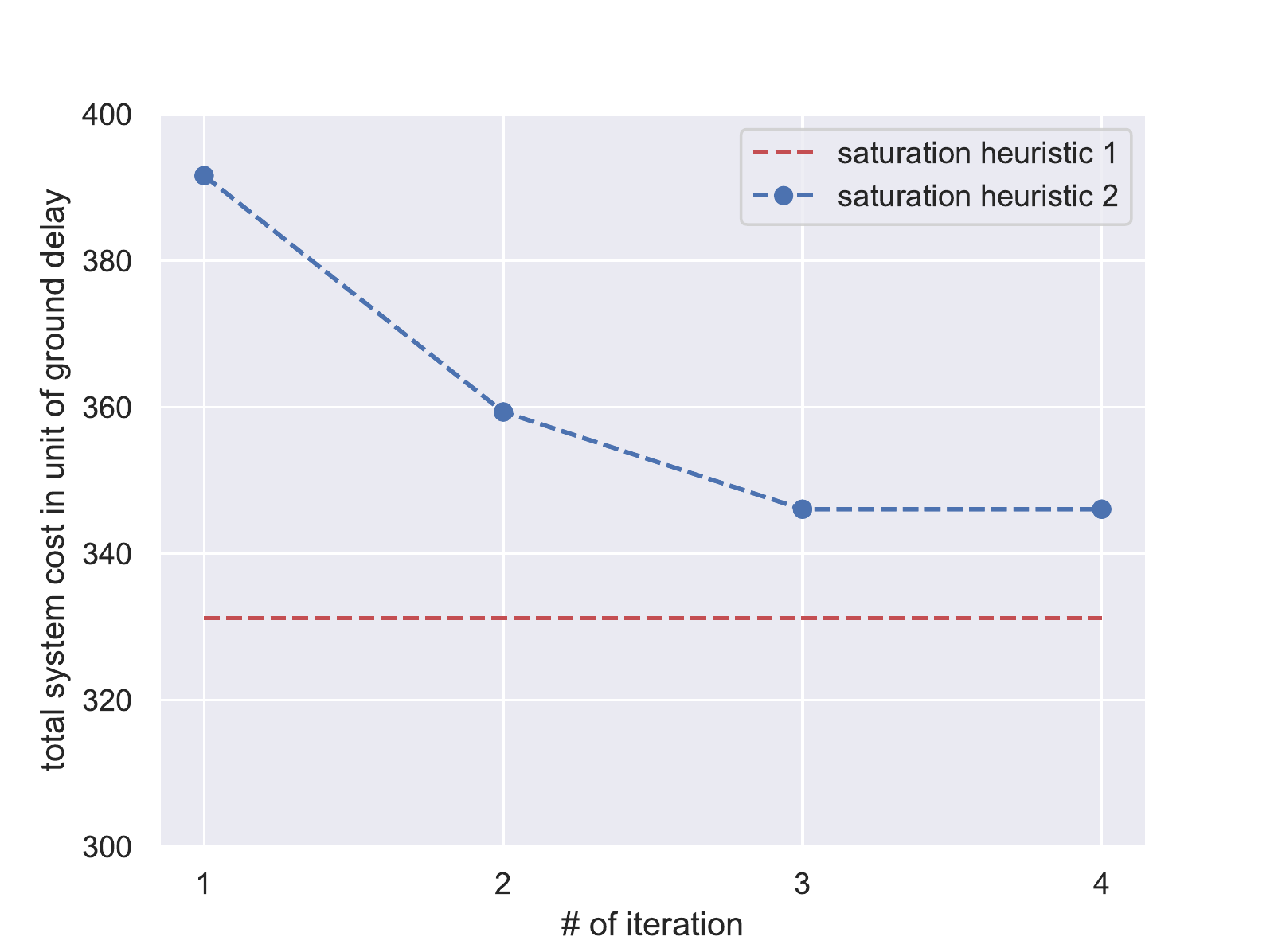}
  \caption{Performance of Saturation Heuristics}\label{Performance_Saturation_Heuristics}
\end{minipage}%
\begin{minipage}{.5\textwidth}
  \centering
  \includegraphics[width=\linewidth]{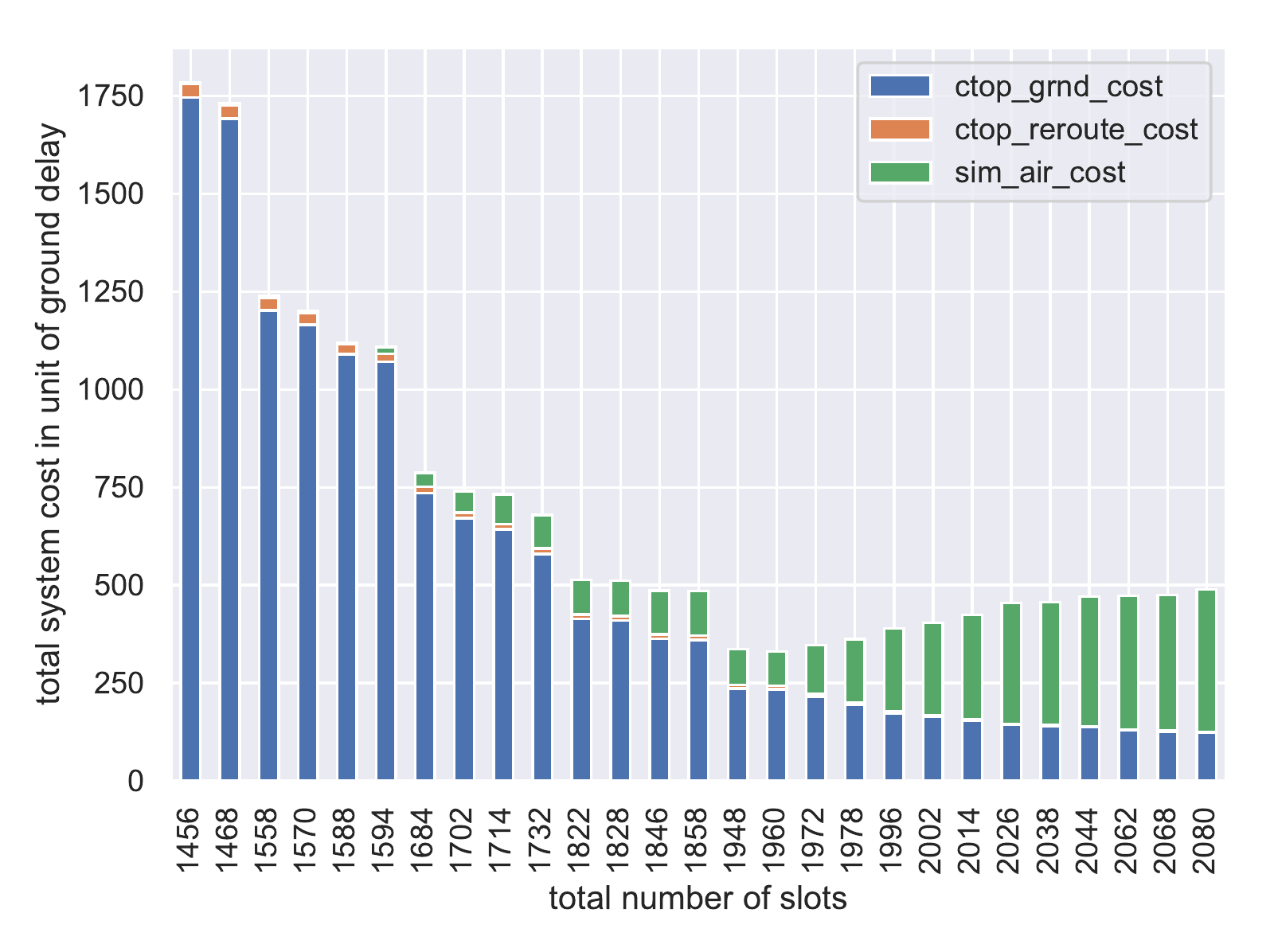}
  \caption{Linear Interpolation of Capacity Scenarios}\label{Linear_Interpolation}
\end{minipage}
\end{figure}

Figure \ref{Performance_Saturation_Heuristics} shows the performance of two phase-1 saturation heuristics described in section \ref{Simulation_based_Optimization}. Heuristic 2 takes four iterations to converge, which performs slightly worst than first saturation heuristic (346.06 versus 331.19). In this use case, the capacity in scenario 1 is strictly better than scenario 2, which is in turn better than scenario 1. Apart from these two saturation heuristics, we also test a third, a more straightforward heuristic: the FCA rates are obtained by doing linear interpolation of scenario capabilities. In Figure \ref{Linear_Interpolation}, the leftmost bar corresponds scenario 3. A lot of ground delay costs are occurred because of the very conservative FCA rate policy. The rightmost bar corresponds to scenario 1. We see a lot more air delay because of the aggressively large FCA rates. It happens in this case the FCA rates obtained from heuristic 1 are the same as the rates from linear interpolation results, which also happen to be the same as the capacity in scenario 2. After we obtain good starting points, in phase 2 we use pattern search to find even better parameters for  EWR and FCA0, because they are the most congested two PCAs.
\begin{itemize}
  \item From the starting point obtained in heuristic 1 and linear interpolation, pattern search converges in 3.7 minutes and reduces the system cost from 331.19 to 324.39 (2.1$\%$ decrease).
  \item From the starting point obtained in heuristic 2, pattern search converges in 3.7 minutes and reduces the system cost from 346.06 to \textbf{324.36} (6.3$\%$ decrease). The solution is shown in Figure \ref{FCA_rates_final_solution}.
\end{itemize}

It can be seen that pattern search method is quite efficient and effective in further reducing system cost, and multiple starting points can help to find the best final solution.

We want to compare the final objective value with benchmarks we have. Two-stage PCA model without rerouting \cite{ZhuAggregate} and two-stage PCA model with rerouting \cite{ZhuCentralized} are models that minimize the system efficiency without considering equity issue. On the other hand, equity is built in CTOP slot allocation algorithm.
\begin{enumerate}
  \item Compared with two-stage PCA model without rerouting whose objective value is 404.0, our model which considers the equity issue achieves a lower system cost. This shows the benefit of allowing rerouting in the face of congestion.
  \item Compared with two-stage PCA model with rerouting whose objective value is 167.07, the objective value of our two-phase optimization framework is almost two times larger. This shows the price of fairness.
\end{enumerate}

\begin{table*}[t]
  \resizebox{1.0\textwidth}{!}{%
  \centering
  \begin{tabular}{cc >{\columncolor[gray]{0.8}}c >{\columncolor[gray]{0.8}}c >{\columncolor[gray]{0.8}} c >{\columncolor[gray]{0.8}}c cccc
  >{\columncolor[gray]{0.8}}c >{\columncolor[gray]{0.8}}c >{\columncolor[gray]{0.8}}c >{\columncolor[gray]{0.8}}c
  cccc
  >{\columncolor[gray]{0.8}}c  >{\columncolor[gray]{0.8}}c >{\columncolor[gray]{0.8}} c >{\columncolor[gray]{0.8}}c c ccc }
  \toprule
  & Resource & 20:00 & 15 & 30 & 45 & 21:00 & 15 & 30 & 45 & 22:00 & 15& 30 & 45 & 23:00 & 15 & 30 & 45 & 00:00 & 15 & 30 & 45 &01:00&15&30&45\\
  \midrule
                         & FCA0 & 13 & 14 & 13 & 13 & 13 & 13 & 13 & 13 & 13 & 13 & 24 & 25 & 25 & 25 & 26 & 25	& 25 & 25 & 25 & 25 & 25& 25 & 25 & 25\\
                         & FCA1 & 44 & 44 & 44 & 44 & 44 & 44 & 44 & 44 & 44 & 44 & 50 & 50 & 50 & 50 & 50 & 50 & 50 & 50 & 50 & 50 & 50 & 50 & 50 & 50\\
                         & FCA2 &  2 &  5 &  5 &  5 &  5 &  5 &  5 &  5 &  5 &  5 &  5 &  5 &  5 &  5 &  5 &  5 &  5 &  5 &  5 &  5 & 5&  5 &  5 &  5 \\
                         & EWR  &  8 &  8 &  9 &  9 &  9 &  7 &  8 &  8 &  8 &  8 & 10 & 11 & 10 & 10 &  9 & 11 & 10 & 10 & 10 & 10 & 10& 10 & 10 & 10 \\
  \bottomrule
  \end{tabular}}
  \caption{FCA Rates Final Solution}\label{FCA_rates_final_solution}
\end{table*}

\section{Conclusions}\label{Conclusions}
The overarching goal of this paper is develop decision making uncertainty algorithms for setting rates across multiple FCAs in CTOP program in the CDM paradigm. To achieve this goal, CTOP FCA rate planning problem has been split into two steps: rate planning given demand estimation, which is a relative easy step and experiences can be borrowed from GDP planning, and rate planning when flight rerouting is also considered, which turns out to be a much more difficult problem.

We first reviewed existing CTOP related stochastic programming models and pointed out the features each model. We believe this categorization of models is enlightening for air traffic flow management research. Second, we focused on FCA-PCA models, which were considered to be promising as they are designed to directly optimize FCA rates given demand information as input. We revealed some of the problems this class of models have due to flow approximation, and showed how these deficiencies can be addressed by PCA models with correctly placed FCAs. Up to this point, we can say that the first step of the problem has been properly solved. However, when we plan to use the first step result in the second step problem, two issues prevent us to achieve optimal FCA rates in the demand uncertainty case, or even obtaining acceptable rate without without using heuristic technique. Third, we proposed a two-phase optimization framework, and combined stochastic optimization with simulation-based optimization. Through a representative use case, we have demonstrated that this framework is efficient and effective.

This is the first fully CDM compatible paper that addresses multiple constrained resources air traffic flow rate planning with reroute options. This work is not only meaningful in providing much-needed decision support capabilities for effective application of CTOP, but also can be .
\section{Acknowledgement}
This research was partly funded under NASA Research Announcement contract \#NNA16BD96C.

\bibliographystyle{elsarticle-harv}
\bibliography{ESOM_bib}

\section{Appendix}
\subsection{Lightweight CTOP Algorithm}
\subsubsection{Assumptions}
\begin{itemize}
  \item There is only one activate TMI, which is the CTOP being implemented. TMI interaction is not considered
  \item Pop-up flights or cancelled flights are not considered
\end{itemize}
\subsubsection{Algorithm Input}
\begin{itemize}
  \item For each FCA $r$, its activation time, acceptance rate $P^r_t$ at each 15 minutes time period and filters 
  \item For each flight, the unimpeded FCA arrival time for each of its TOS option
\end{itemize}
\subsubsection{Slot Creation}
We create evenly spaced slots for each time period based on given acceptance rate. The $i$-th slot in time period $t$ FCA $r$ is:
\begin{equation}
\begin{split}
& \text{slot}^r_t(i) = \text{round} ((i-1)\times \frac{15\times 60}{P^r_t}) 
\end{split}
\end{equation}

\subsubsection{Slot Allocation Algorithm}\label{slot_allocation_algorithm}
{
\begin{algorithm} \caption{CTOP TOS Allocation Algorithm \cite{OptimalAirlineActions}}
\begin{algorithmic}[1]
\STATE Determine flights included by the CTOP program. A flight is included in CTOP if any TOS route intersect any of the CTOP's FCAs during active periods
\STATE Determine flights that are part of CTOP demand but are exempted
\STATE Assign slots to exempted flights first
\STATE Sort flights by Initial Arrival Time (IAT), which is the earliest FCA arrival time at any of a CTOP's FCAs using any of the flight's TOS options
\STATE Once at a time, in IAT order, assign each flight the lowest adjusted cost trajectory and slot
\end{algorithmic}
\end{algorithm}
}

If a TOS route intersect two or more FCAs. Take the second FCA as an example: slot will be marked as used by finding the first available slot in this FCA that has a time equal or later than the time at the flight would intersect this FCA if flight departing at its ETD, which would include any delay first (primary) FCA imposes.

\subsubsection{Stochastic Flow Simulation Algorithm}\label{flow_simulation_algo}
We will only consider the cost for non-exempted flights. The total cost is composed of three parts: reroute cost, ground delay cost incurred by FCAs, air delay cost incurred by PCAs capacity constraints. Like in all FCA-PCA models, we require that FCAs must be placed before PCAs. After we run the CTOP allocation algorithm, we can easily calculate the costs of first two parts, and we will know $S^k_{t,\rho}$. By solving the following optimization problem, we will know the third part.
\begin{align}
& \min\quad c_a\sum_{q\in \mathcal{Q}}p_q\sum_{t\in \mathcal{T}}\sum_{\rho \in \mathscr{P}}\sum_{k \in \rho}A^{k,q}_{t,\rho} \\
& L^{k,q}_{t,\rho}            = \begin{cases}
                             \mbox{if $k=\rho_1$} &   S^k_{t,\rho}  -(A^{k,q}_{t,\rho}-A^{k,q}_{t-1,\rho})\\
                             \mbox{else}       &   \text{UpPCA}^{k,q}_{t,\rho}-(A^{k,q}_{t,\rho}-A^{k,q}_{t-1,\rho})\\
                           \end{cases} \label{TwoStagePass} & \forall t\in \mathcal{T}, q\in \mathcal{Q}, \rho\in \mathscr{P}, k\in \rho\\
& \text{UpPCA}^{k,q}_{t,\rho} = L^{{k^\prime},q}_{t-\Delta^{k^\prime,k},\rho} & \forall t\in \mathcal{T}, q\in \mathcal{Q}, (k^\prime,k)\in \rho \label{TwoStageUpPCA} \\
& \sum_{\rho \in \mathscr{P}}L^{k,q}_{t,\rho} \le M^k_{t,q} &  \forall t\in \mathcal{T},q\in \mathcal{Q}, k\in \mathcal{P}\\
& L^{k,q}_{t,\rho}, A^{k,q}_{t,\rho} \in \mathbb{Z}_+ & \forall t\in \mathcal{T},q\in \mathcal{Q},\rho\in \mathscr{P},k\in \rho
\end{align}
\section{Appendix for Acronyms and Abbreviations}
\noindent
  \textbf{TMI}   \hfill{Traffic Management Initiative}\\
  \textbf{GDP}   \hfill{Ground Delay Program}\\
  \textbf{AFP}   \hfill{Airspace Flow Program}\\
  \textbf{FCA}   \hfill{Flow Constrained Area}\\
  \textbf{PCA}   \hfill{Potential Constrained Area}\\
  \textbf{CDM}   \hfill{Collaborative Decision Making}\\
  \textbf{CTOP}  \hfill{Collaborative Trajectory Options Program}\\
  \textbf{TOS}   \hfill{Trajectory Options Set}\\
  \textbf{ESOM}  \hfill{Enhanced Stochastic Optimization Model}
\end{document}